\documentclass[a4paper,10pt,reqno]{amsart}
\usepackage{amsmath}
\usepackage{cases}

\usepackage{stmaryrd}

\usepackage{amsfonts}
\usepackage[colorlinks,linkcolor=blue,citecolor=blue]{hyperref}
\usepackage{latexsym, amssymb, amsmath, amsthm, bbm}
\usepackage[all]{xy}
\usepackage{pgfplots}
\usepackage{enumerate}

\DeclareSymbolFont{EulerExtension}{U}{euex}{m}{n}
\DeclareMathSymbol{\euintop}{\mathop} {EulerExtension}{"52}
\DeclareMathSymbol{\euointop}{\mathop} {EulerExtension}{"48}

\allowdisplaybreaks[4]

\def \cop{\mathrm{cop}}

\def \id{\operatorname{id}}

\def \C{\mathcal{C}}

\def \M{\mathrm{M}}
\def \N{\mathbb{N}}

\def \k{\Bbbk}

\def \dim{\operatorname{dim}}

\def \C{\mathcal{C}}
\def \D{\Delta}

\def \M{\mathrm{M}}
\def \N{\mathbb{N}}

\def \cop{\mathrm{cop}}

\def \A{\mathcal{A}}
\def \B{\mathcal{B}}
\def \C{\mathcal{C}}
\def \D{\mathcal{D}}
\def \E{\mathcal{E}}

\def \X{\mathcal{X}}
\def \Y{\mathcal{Y}}

\def \W{\mathcal{W}}
\def \M{\mathcal{M}}
\def \N{\mathcal{N}}

\usepackage{mathtools}

\numberwithin{equation}{section}

\newtheorem{theorem}{Theorem}[section]
\newtheorem{lemma}[theorem]{Lemma}
\newtheorem{proposition}[theorem]{Proposition}
\newtheorem{corollary}[theorem]{Corollary}
\newtheorem{definition}[theorem]{Definition}

\newtheorem{remark}[theorem]{Remark}

\newtheorem{notation}[theorem]{Notation}%[section]

\begin{document}
\title[A matrix decomposition for coalgebras with applications]{A matrix decomposition for coalgebras with applications to antipode formulas}

\author[J.-W. He]{Ji-Wei He}
\address{School of Mathematics, Hangzhou Normal University, Hangzhou 311121, China}
\email{jwhe@hznu.edu.cn}

\author[X. Kong]{Xiaojie Kong}
\address{School of Mathematics, Hangzhou Normal University, Hangzhou 311121, China}
\email{xjkong@stu.hznu.edu.cn}

\author[K. Li]{Kangqiao Li${}^\dag$}
\address{School of Mathematics, Hangzhou Normal University, Hangzhou 311121, China}
\email{kqli@hznu.edu.cn}



\thanks{$\dag$ Corresponding author}

\thanks{The work was partially supported by NSFC [Nos. 12301049, 12371017].}

\date{}

\begin{abstract}
In this paper, we provide a matrix basis decomposition for $\mathbb{N}$-graded coalgebras over an algebraically closed field,
and establish a formula on the coproducts of the matrix basis.
Our main result can be applied to construct linear operators from such coalgebras.
As an example,
%As an application, we consider weak bialgebras whose underlying coalgebras admit such a grading.
we obtain recursive formulas extending the antipode of the coradical of a weak bialgebra, thereby making it a weak Hopf algebra.
%
%that extend a linear map on the degree-zero component to an antipode on the whole weak bialgebra, thereby making it a weak Hopf algebra.
\end{abstract}

\maketitle

\section{Introduction}
%Primitive matrices were used by Li and Zhu in their study of exponents of non-cosemisimple Hopf algebras \cite{LZ19}. Given basic multiplicative matrices $\C$ and $\D$, a $(\C, \D)$-primitive matrix $\W$ satisfies
%$$
%\Delta(\W)=\C\;\widetilde{\otimes}\;\W
%+\W\;\widetilde{\otimes}\;\D
%$$
%where $\;\widetilde{\otimes}\;$ denotes the matrix tensor product: the usual rule for matrix multiplication is used, with each product of entries replaced by their tensor product. This identity extends the familiar coproduct formula for skew-primitive elements to the setting of arbitrary simple subcoalgebras. See  \cite[Definition 3.2]{LZ19}. The relationship between primitive matrices and simple bicomodules was further developed by Yu, Li and Liu \cite{YLL24}, who used complete families of non-trivial primitive matrices to describe $H_1/H_0$ and study link quivers of Hopf algebras with the dual Chevalley property. Here the entries of the primitive matrices represent basis elements in the corresponding quotients.
Over an algebraically closed field, every simple coalgebra admits a basis through which the comultiplication can be written in a matrix form. Basic multiplicative matrices therefore provide a convenient way to describe simple subcoalgebras and coactions of simple comodules. Multiplicative matrices were introduced in \cite{Man88}, and primitive matrices were subsequently considered in \cite{LZ19}.

For a primitive matrix, the comultiplication consists of two terms determined by a pair of basic multiplicative matrices. In higher homogeneous degrees, additional terms involving positive lower degrees must be taken into account. This leads to the question of how to choose bases in matrix form so that comultiplication can be described uniformly in every degree.

In this paper, we study this questions for $\mathbb{N}$-graded coalgebras over an algebraically closed field $\k$. More precisely, we assume that
$$
H=\bigoplus_{n\geq0}H(n)\;\;\;\;\text{and}\;\;\;\;H(0)=H_0,
$$
where $H_0$ is the coradical of $H$. Each homogeneous component $H(n)$ is naturally an $H(0)$-bicomodule and is a direct sum of simple bicomodules. For simple coalgebras $C$ and $D$ with basic multiplicative matrices $\C$ and $\D$, respectively, we describe a basis of a simple $C$-$D$-bicomodule by a matrix $\M$ satisfying
$$
\rho^L(\M)=\C\;\widetilde{\otimes}\;\M\;\;\;\;\text{and}\;\;\;\;
\rho^R(\M)=\M\;\widetilde{\otimes}\;\D,
$$
where the matrix tensor product is the one defined in Notation \ref{notation}. We also describe bases of cotensor products by matrix tensor products. Combining these results gives the following matrix basis decomposition.

\begin{theorem}
Fix a basic multiplicative matrix for each simple subcoalgebra of $H$, and let $\mathcal{P}_0$ be the set of all such matrices. Then for each $n\geq1$, there exists a set $\mathcal{P}_{n}$ of matrices over $H(n)$ satisfying the following properties:
\begin{itemize}
\item[(1)]  Every matrix in $\mathcal{P}_n$ has linearly independent entries.
\item [(2)] For any $\C, \D\in \mathcal{P}_0$, let $C=\mathrm{span}(\C)$ and $D=\mathrm{span}(\D)$. There exists a subset ${}^\C{\mathcal{P}_n}{}^\D\subseteq\mathcal{P}_n$ such that
    $$
    {}^C{H(n)}{}^D=\bigoplus_{\mathcal{W}
    \in{}^\C{\mathcal{P}_n}^\D}\mathrm{span}(\mathcal{W}).
    $$
\item [(3)] $\mathcal{P}_n$ is the disjoint union:
    $$
    \mathcal{P}_n=\bigcup\limits_{\C, \D\in\mathcal{P}_0}
    {}^\C{\mathcal{P}_n}{}^\D.
    $$
\item [(4)] For any $\C, \D\in \mathcal{P}_0$ and every $\mathcal{W}\in {}^\C{\mathcal{P}_n}^\D$, there exist $\alpha_{\mathcal{X},\mathcal{Y}}\in\k$, such that
    \begin{eqnarray}\label{1}
    \Delta(\mathcal{W})=\C\;\widetilde{\otimes}\;\mathcal{W}
    +\mathcal{W}\;\widetilde{\otimes}\;\D
    +\sum_{i=1}^{n-1}\sum_{\mathcal{E}\in\mathcal{P}_0}
    \sum_{\substack{\mathcal{X}\in{}^\C{\mathcal{P}_i}^\mathcal{E},\\
    \mathcal{Y}\in{}^\mathcal{E}{\mathcal{P}_{n-i}}^\D}}
    \alpha_{\mathcal{X},\mathcal{Y}}(\mathcal{X}\;\widetilde{\otimes}\;
    \mathcal{Y}),
    \end{eqnarray}
where the coefficients $\alpha_{\mathcal{X},\mathcal{Y}}$ are uniquely determined.
\end{itemize}
\end{theorem}
The theorem extends this description to all homogeneous degrees. Coassociativity gives a further identity for the lower-degree terms, as stated in Corollary \ref{coassW}.

We apply the matrix basis decomposition to weak bialgebras. Weak Hopf algebras relax the usual compatibility conditions involving the unit and counit; see \cite{BNS99} and \cite{S03}. Their antipodes satisfy three convolution identities involving the source and target counital maps. Thus a recursive construction must preserve all three identities at each step. We aim to obtain such a construction directly from the matrix formulas for comultiplication.

The main application is the following extension theorem.
We remark that it is a consequence of a talk once given by H. Huang and her collaborators X. Wang and J.J. Zhang, while our construction with recursive formulas may be helpful in calculations.

\begin{theorem}\label{them}
Let $H=\bigoplus_{n\geq0}H(n)$ be a weak bialgebra with an $\mathbb{N}$-grading as coalgebra, and let $H(0)=H_0$.
Suppose that there is a linear map
$S: H(0)\rightarrow H$ such that
$$
S\ast\id=\varepsilon_s,\;\;\id\ast S=\varepsilon_t,\;\;
S\ast\id\ast S=S,
$$
where $\varepsilon_s$ and $\varepsilon_t$ are the source and target counital maps, respectively.
Then $H$ is a weak Hopf algebra.
\end{theorem}
The proof is constructive. Suppose that $S$ has already been defined on the homogeneous components of degrees less than $n$ and satisfies the three antipode identities in Theorem \ref{them}. Applying (\ref{1}), we obtain two equivalent expressions for the values of $S$ on the entries of a basis matrix in degree $n$. Repeating this construction gives an antipode on the entire weak bialgebra.

The paper is organized as follows. In Section 2, we recall the notations of bicomodules, cotensor products, graded coalgebras, and multiplicative matrices. In Section 3, we first give matrix descriptions of bases of simple bicomudules and their cotensor products. Then we prove the matrix basis decomposition and derive the identity arising from coassociativity. In Section 4, the main results are applied to bialgebras. Then we obtain the recursive formulas for the antipode, and prove the extension theorem.

\section{Preliminaries}\label{section2}

We refer to \cite{Swe69, Mon93, Rad12}  for the definitions and basic properties about coalgebras.

Throughout this paper, all vector spaces and coalgebras are defined over an algebraically closed field $\k$, and the tensor product over $\k$ is denoted by $\otimes$ for simplicity.
For a vector space $V$, we always denote by $\mathbf{M}_{r\times s}(V)$ the space of $r\times s$ matrices over $V$, and write $\mathbf{M}_{r\times r}(V)$ as $\mathbf{M}_r(V)$.

For vector spaces $U$ and $V$, let $\tau_{U, V}: U\otimes V\rightarrow V\otimes U$ denote the flip map.
For a coalgebra $C$,  we denote its coopposite coalgebra by $C^{\cop}$ with comultiplication $\Delta_C^{\cop}=\tau_{C, C}\circ\Delta_C$.

\subsection{Bicomodules and cotensor products}

We begin by recalling the notion of bicomodules over coalgebras.

\begin{definition}
Let $C$ and $D$ be coalgebras. A $C$-$D$-bicomodule is a vector space $M$ equipped with a left $C$-comodule structure $\rho^L: M\rightarrow C\otimes M$ and a right $D$-comodule structure $\rho^R: M\rightarrow M\otimes D$ satisfying
$$
(\rho^L\otimes\id_D)\circ\rho^R=(\id_C\otimes\rho^R)\circ\rho^L.
$$
\end{definition}

We denote by ${}^C\mathfrak{M}{}^D$ the $\k$-linear abelian category of $C$-$D$-bicomodules.
Moreover, the category ${}^C\mathfrak{M}{}^D$ can be identified with the category ${}^{C\otimes D^{\cop}}\mathfrak{M}$ of left $C\otimes D^{\cop}$-comodules:

%The following lemma provides a useful identification between bicomodules and comodules over a tensor product coalgebra.

\begin{lemma}\label{cong}
Let $C$ and $D$ be coalgebras. Then there is an isomorphism of $\k$-linear abelian categories
$$
\Phi: {}^C\mathfrak{M}^D\cong {}^{C\otimes D^\cop}\mathfrak{M},\;\;\; M\mapsto M,
$$
which sends a $C$-$D$-bicomodule $M$ to the same vector space with a left $C\otimes D^{\cop}$-coaction given by
$$
\delta: M\rightarrow (C\otimes D^\cop)\otimes M,\;\;\;\; m\mapsto\sum(m_{(-1)}\otimes m_{(1)})\otimes m_{(0)}.
$$
Here, we use Sweedler notation $m\mapsto\sum m_{(-1)}\otimes m_{(0)}$ for the left $C$-coaction and $m\mapsto\sum m_{(0)}\otimes m_{(1)}$ for the right $D$-coaction.
\end{lemma}

\begin{proof}
It is straightforward to verify that $(M, \delta)$ is a left $C\otimes D^{\cop}$-comodule and $\Phi$ is a $\k$-linear (addtive) functor.

Moreover, the inverse functor $\Psi: {}^{C\otimes D^\cop}\mathfrak{M}\rightarrow{}^C\mathfrak{M}^D$ of $\Phi$ is defined as follows.
%The inverse functor $\Psi$ sends a left $C\otimes D^\cop$-comodule $(M, \mu)$ to the $(C, D)$-bicomodule
For each object $(M, \delta)\in {}^{C\otimes D^\cop}\mathfrak{M}$, the left $C$-coaction $\rho^L$ and right $D$-coaction $\rho^R$ on $\Psi(M)$ are defined by
\begin{equation}\label{rhoL}
\rho^L(m)=(\id_C\otimes\varepsilon_D\otimes\id_M)\circ\delta(m)
\in C\otimes M,
\end{equation}
\begin{equation}\label{rhoR}
\rho^R(m)=\tau_{D,M}\circ(\varepsilon_C\otimes\id_D\otimes\id_M)\circ\delta(m)
\in M\otimes D
\end{equation}
for all $m\in M$.
\end{proof}

%Denote by $\mathcal{S}$ the set of all the simple subcoalgebras of a coalgebra $H$.
By this lemma, one knows that if $C$ and $D$ are cosemisimple coalgebras, then ${}^C\mathfrak{M}{}^D$ is semisimple.

In particular, we obtain the following result on simple bicomodules over simple coalgebras.

\begin{corollary}\label{dim=rs}
Let $C$ and $D$ be simple coalgebras with $\dim(C)=r^2$ and $\dim(D)=s^2$. Then
%\begin{itemize}
%\item[\textup{(1)}](\cite[Lemma 2.3]{YLL24})\label{YLL}
%Every simple $C$-$D$-bicomodule has dimension $rs$.
%\item[\textup{(2)}]
%Every $C$-$D$-bicomodule is a direct sum of simple bicomodules. In particular,
a $C$-$D$-bicomodule $M$ is simple if and only if $\dim(M)=rs$.
\end{corollary}

%\begin{proof}
%Since $C\otimes D^{\cop}$ is also a simple coalgebra, the desired claim follows from Lemma \ref{cong}.
%\end{proof}

Now we recall the notion of the cotensor product of comodules.

\begin{definition}(\cite[Section 0]{T77})
Let C be a coalgebra. For a right $C$-comodule $M$ and a left $C$-comodule $N$, their cotensor product $M\square_C N$ is the kernel of
$$
\rho_M\otimes\id-\id\otimes\rho_N: M\otimes N\rightarrow M\otimes C\otimes N,
$$
where $\rho_M$ and $\rho_N$ denote the structure maps of $M$ and $N$, respectively.
\end{definition}

%Let $\mathcal{C}$ denote a symmetric monoidal category and let $M\in\mathcal{C}$. If $M$ is a left $H$-module and right $H$-module,

Furthermore, for a finite-dimensional right $C$-comodule $(M, \rho_M)$ and a finite-dimensional left $C$-comodule $(N, \rho_N)$, there is a natural isomorphism
\begin{equation}\label{MN=MastNast}
(M\square_C N)^\ast\cong M^\ast\otimes_{C^\ast} N^\ast
\end{equation}
according to \cite[Lemma 3]{CKQ05}, where $(M^\ast, \leftharpoonup)$ and $(N^\ast, \rightharpoonup)$ are the induced right and left $C^\ast$-module structures, respectively. Specifically, for any $m\in M$, $n\in N$, if we write
$$\rho_M(m)=\sum m_{(0)}\otimes m_{(1)}\;\;\;\;
\text{and}
\;\;\;\;\rho_N(n)=\sum n_{(-1)}\otimes n_{(0)},
$$
then these module structures are given by
$$
(\varphi\rightharpoonup f)(m)=\sum\varphi(m_{(0)})f(m_{(1)})\;\;\;\;\text{and}\;\;\;\;
(f\leftharpoonup\psi)(n)=\sum f(n_{(-1)})\psi(n_{(0)}),
$$
for all $f\in C^\ast$, $\varphi\in M^\ast$, $\psi\in N^\ast$.

We use the isomorphism $(\ref{MN=MastNast})$ to prove the following lemma.

\begin{lemma}\label{cotenp dim}
Let $C$ be a simple coalgebra. Suppose $M$ is a finite-dimensional right $C$-comodule and let $N$ be a finite-dimensional left $C$-comodule. Then
$$
\dim(M\square_C N)=\frac{\dim(M)\dim(N)}{\dim(C)}.
$$
\end{lemma}

\begin{proof}
Note that $C$ is a simple coalgebra over an algebraically closed field $\k$.
Thus $C^\ast\cong \mathbf{M}_r(\k)$ as algebras for some $r\geq1$. It is clear that
$\mathbf{M}_{1\times r}(\k)$ and $\mathbf{M}_{r\times 1}(\k)$ are respectively the simple right and left $\mathbf{M}_r(\k)$-modules, and then there exist positive integers $k$ and $l$ such that
$$
M^\ast\cong \mathbf{M}_{1\times r}(\k)^{\oplus k},\;\;\;\;\text{and}\;\;\;\;N^\ast\cong \mathbf{M}_{r\times 1}(\k)^{\oplus l},
$$
as $C^\ast$-modules.
It follows that $M^\ast\otimes_{C^\ast}N^\ast\cong
\big(\mathbf{M}_{1\times r}(\k)\otimes_{\mathbf{M}_r(\k)}\mathbf{M}_{r\times 1}(\k)\big)^{\oplus kl}$.

Moreover, one can directly find that $\mathbf{M}_{1\times r}(\k)\otimes_{\mathbf{M}_r(\k)}\mathbf{M}_{r\times 1}(\k)$
is 1-dimensional.
Consequently, we know by (\ref{MN=MastNast}) that
\begin{eqnarray*}
\dim(M\square_CN)=\dim(M^\ast\otimes_{C^\ast}N^\ast)=kl
=\frac{kr\cdot lr}{r^2}
=\frac{\dim(M)\dim(N)}{\dim(C)}.
\end{eqnarray*}

\end{proof}

\subsection{$\mathbb{N}$-graded coalgebras}

In this subsection, we recall the definition of coradically graded coalgebras and the corresponding coradical orthonormal idempotents in its dual algebras.
\begin{definition}\label{gcoalg}
An $\mathbb{N}$-graded coalgebra is a coalgebra $H$ with a decomposition $H=\bigoplus_{n\geq0}H(n)$
such that
$$
\Delta(H(n))\subseteq\bigoplus\limits_{i=0}^nH(i)\otimes H(n-i)\;\;\;\;(n\geq0),\;\;\;\;\text{and}\;\;\;\;
\varepsilon(H(n))=0\;\;\;\;(n>0).
$$
\end{definition}

\begin{remark}
In the situation of Definition \ref{gcoalg}, the coalgebra $H(0)$ is cosemisimple if and only if $H(0)$ is the coradical $H_0$ of $H$.
\end{remark}

%Write $H_0$ for the coradical of $H$, namely the sum of all simple subcoalgebras of $H$, and let $(H_n)_{n\geq0}$ be the coradical filtration, with
%$$
%H_n=\Delta^{-1}(H_0\otimes H+H\otimes H_{n-1})\;\;\;\;(n\geq1).
%$$
%The coalgebra is called coradically graded if $H_n=\bigoplus\limits_{j=0}^n H(j)$ for all $n\geq0$. In particular, $H(0)=H_0$.

For the remainder of this paper, we always assume that $H$ is a coalgebra, and denote by $\mathcal{S}$ the set of all the simple subcoalgebras of $H$.

\begin{notation}\label{eC}
Let $H=\bigoplus_{n\geq0}H(n)$ be an $\mathbb{N}$-graded coalgebra with $H(0)=H_0$.

\begin{itemize}
\item[\textup{(1)}]
For any $C\in\mathcal{S}$, define a linear function $e_C\in H^\ast$ as follows:
$$e_C|_D=\delta_{C,D}\varepsilon|_D \;\;\;\;(\forall D\in\mathcal{S})\;\;\;\;\text{and}\;\;\;\;e_C|
_{\bigoplus_{n\geq1}H(n)}=0.
$$
\item[\textup{(2)}]
For any $h\in H$ and $C, D\in\mathcal{S}$, we write
$${^C}h=h\leftharpoonup e_C,\;\;\;h^D=e_D\rightharpoonup h,\;\;\;{^C}h^D=e_D\rightharpoonup h\leftharpoonup e_C,$$
where $\leftharpoonup$ and $\rightharpoonup$ are hit actions of $H^\ast$ on $H$. Notations such as $V^{C}:=e_C\rightharpoonup V$ for a subspace $V$ of $H$ are used as well.
\end{itemize}
\end{notation}

\begin{remark}
It is straightforward to find that
$$e_Ce_D=\delta_{C,D}e_C\;\;\;\;
(\forall\;C,D\in\mathcal{S})\;\;\;\;\text{and}\;\;\;\;\sum_{C\in\mathcal{S}}e_C=\varepsilon.$$
hold, where the sum $\sum_{C\in\mathcal{S}}e_C$ is well-defined in $H^\ast$ (see \cite[Section 1]{Rad78}). Thus $\{e_C\}_{C\in\mathcal{S}}$ is referred to as a family of coradical orthonormal idempotents by \cite{LZ19}.
\end{remark}

Some direct properties, which are analogous to \cite[Proposition 2.6]{LL22}, are listed as follows:
%\begin{proposition}
%Let $H$ be a coalgebra. Then for all $C, D\in\mathcal{S}$, we have
%\begin{itemize}
%\item[(1)] ${}^CH_0{}^D=\delta_{C, D}C$;
%\item[(2)] ${}^CH_1{}^D\subseteq\Delta^{-1}(C\otimes {}^CH_1{}^D+{}^CH_1{}^D\otimes D)$;
%\item[(3)] ${}^CH{}^D\subseteq \Ker(\varepsilon)$ if $C\neq D$;
%\item[(4)] Suppose $V\subseteq H$ is a $\k$-subspace. The following direct-sum decompositions hold:
%    \begin{itemize}
%    \item[(i)] $V=\bigoplus\limits_{C\in\mathcal{S}}{}^CV$ if $V$ is a left coideal;
%    \item[(ii)] $V=\bigoplus\limits_{D\in\mathcal{S}}V{}^D$ if $V$ is a right coideal;
%    \item[(iii)] $V=\bigoplus\limits_{C,D\in\mathcal{S}}{}^CV{}^D$ if $V$ is a subcoalgebra.
%    \end{itemize}
%\end{itemize}
%\end{proposition}

\begin{proposition}\label{Hn}
Let $H=\bigoplus_{n\geq0}H(n)$ be an $\mathbb{N}$-graded coalgebra with $H(0)=H_0$, and let $\{e_C\}_{C\in\mathcal{S}}$ be the family defined in Notation \ref{eC}. For $C, D\in\mathcal{S}$ and $n\geq0$, write
$$
{}^CH(n){}^D=e_D\rightharpoonup H(n)\leftharpoonup c_C.
$$
Then we have
\begin{itemize}
\item[(1)] ${}^CH(0){}^D=\delta_{C, D}C.$
\item[(2)] For every $n\geq0$,
$$
\Delta({}^CH(n){}^D)\subseteq\bigoplus\limits_{i=0}^n
\bigoplus\limits_{E\in\mathcal{S}}{}^CH(i){}^E\otimes
{}^EH(n-i){}^D.
$$
%In particular,
%$$\Delta({}^CH(1){}^D)\subseteq C\otimes{}^CH(1){}^D+{}^CH(1){}^D\otimes D.
%$$
\item[(3)]
$\varepsilon({}^CH(n){}^D)=0$ whenever $n\geq1$ or $C\neq D$.
\item[(4)]
$
H(n)=\bigoplus_{C, D\in\mathcal{S}}{}^CH(n){}^D.
$
\end{itemize}
\end{proposition}

\subsection{Multiplicative matrices over coalgebras}
The notion of multiplicative matrices over coalgebras was once introduced in \cite{Man88}, and it was used to study non-pointed Hopf algebras. In this subsection, let us recall the definition and collect some basic properties for later use.

\begin{notation}\label{notation}
Let $U$ and $V$ be vector spaces.
\begin{itemize}
\item[(1)]
For matrices $\A:=(u_{ij})_{m\times n}$ over $U$ and $\B:=(v_{ij})_{n\times l}$ over $V$, denote the following matrix
$$\A\;\widetilde{\otimes}\;\B:=\left(\sum_{k=1}^n u_{ik}\otimes v_{kl}\right)_{m\times l};$$
\item[(2)]
For a matrix $\A:=(u_{ij})_{m\times n}$ over $U$, we denote its transpose by
$\A^\mathrm{T}:=(u_{ji})_{n\times m};$
\item[(3)]
For any linear map $f: U\rightarrow V$ and a matrix $\mathcal{A}=(u_{ij})_{m\times n}$ over $U$, denote the following matrix
$$
f(\mathcal{A})=(f(u_{ij}))_{m\times n}.
$$
\end{itemize}
\end{notation}

\begin{remark}\label{AB}
Under the assumptions in Notation \ref{notation}(1), it is straightforward to verify the following facts:
\begin{itemize}
\item[(1)]
 If the entries of $\A$ and $\B$ are linearly independent, respectively, then $\A\;\widetilde{\otimes}\;\B$ has linearly independent entries;
\item[(2)] If $\B$ has linearly independent entries, then $\A\;\widetilde{\otimes}\;\B=0$ implies $\A=0$.
\end{itemize}
\end{remark}

For a positive integer $r$, denote by $I_r$ the $r\times r$ identity matrix over $\k$.

\begin{definition}(\cite{Man88} and \cite[Section 3]{LZ19})\label{MM}
Let $(H,\Delta,\varepsilon)$ be a coalgebra over $\k$.
\begin{itemize}
  \item[(1)] A square matrix $\C_{r\times r}$ over $H$ is called multiplicative, if $\Delta(\C_{r\times r})=\C_{r\times r}\;\widetilde{\otimes}\;\C_{r\times r}$ and $\varepsilon(\C_{r\times r})=I_r$.
      Indeed,  $\Delta$ and $\varepsilon$ act entrywise. Equivalently, its entries satisfy
$$
\Delta(c_{ij})=\sum_{k=1}^r c_{ik}\otimes c_{kj},\;\;\;\;\varepsilon(c_{ij})=\delta_{ij}\;\;\;\;
(1\leq i, j\leq r).
$$
  \item[(2)] A multiplicative matrix is called basic if its entries are linearly independent.
\end{itemize}
\end{definition}

\begin{remark}
Suppose that $\C$ is a basic multiplicative matrix over a coalgebra $H$. Then the entries of $\C$ span a simple subcoalgebra $C$ of $H$. In this case, we say that $\C$ is a basic multiplicative matrix of $C$.
\end{remark}
%
%Two square matrices $\C$ and $\D$ of the same size over a vector space $V$ are similar if there exists an invertible matrix $L$ over $\k$ such that $\D=L\C L^{-1}$. This is denoted by $\C\sim\D$ for simplicity.

\begin{lemma}\label{BMM}
Let $C$ be a simple coalgebra. Then:
\begin{itemize}
\item[(1)]
There exists a basic multiplicative matrix $\C$ of $C$.
\item[(2)](\cite[Lemma 2.4]{L22}) Let $\C$ be a basic multiplicative matrix of $C$. Then $\D$ is also a basic multiplicative matrix of $C$ if and only if $\D=P\C P^{-1}$ for some invertible matrix $P$ over $\k$.
\end{itemize}
\end{lemma}

%The next lemma shows that the transpose of a multiplicative matrix over $C$ is multiplicative over $C^{\cop}$.

\begin{proposition}\label{Delta(MM)}
Let $C$ be a coalgebra and let $\C$ be a multiplicative matrix over a coalgebra $H$.
%such that
%$$
%\Delta(c_{ij})=\sum_{k=1}^r c_{ik}\otimes c_{kj}\;\;\;\;
%\text{and} \;\;\;\;\varepsilon(c_{ij})=\delta_{ij}
%\;\;\;\;(1\leq i, j\leq r).
%$$
Then $\C^\mathrm{T}$ is a multiplicative matrix over $H^\cop$.

\end{proposition}

\begin{proof}
Suppose $\C=(c_{ij})_{r\times r}$, and we can infer by Definition \ref{MM}(1) that
%for $1\leq i,j \leq r$,
%$
%\Delta(c_{ij})=\sum_{k=1}^r c_{ik}\otimes c_{kj}.
%$
%Equivalently, the above identities can be written in matrix form as
%$$
%\Delta(\mathcal{C})=
%\left(\begin{array}{cccc}
%    c_{11} & c_{12} & \cdots & c_{1r}  \\
%    c_{21} & c_{22} & \cdots & c_{2r}  \\
%    \vdots & \vdots & \ddots & \vdots  \\
%    c_{r1} & c_{r2} & \cdots & c_{rr}
%\end{array}\right)
%\;\widetilde{\otimes}\;
%\left(\begin{array}{cccc}
%    c_{11} & c_{12} & \cdots & c_{1r}  \\
%    c_{21} & c_{22} & \cdots & c_{2r}  \\
%    \vdots & \vdots & \ddots & \vdots  \\
%    c_{r1} & c_{r2} & \cdots & c_{rr}
%\end{array}\right)
%=
%\mathcal{\C} \;\widetilde{\otimes}\; \mathcal{\C}.
%$$
$$\Delta^{\cop}(c_{ij})=\sum_{k=1}^r c_{kj}\otimes c_{ik}\;\;\;\; (\forall 1\leq i,j \leq r)
$$
holds on $C^{\cop}$.
In other words,
$$
\Delta^\cop(\mathcal{C^\mathrm{T}})=
\left(\begin{array}{cccc}
    c_{11} & c_{21} & \cdots & c_{r1}  \\
    c_{12} & c_{22} & \cdots & c_{r2}  \\
    \vdots & \vdots & \ddots & \vdots  \\
    c_{1r} & c_{2r} & \cdots & c_{rr}
\end{array}\right)
\;\widetilde{\otimes}\;
\left(\begin{array}{cccc}
    c_{11} & c_{21} & \cdots & c_{r1}  \\
    c_{12} & c_{22} & \cdots & c_{r2}  \\
    \vdots & \vdots & \ddots & \vdots  \\
    c_{1r} & c_{2r} & \cdots & c_{rr}
\end{array}\right)
=
\mathcal{\C^\mathrm{T}} \;\widetilde{\otimes}\; \mathcal{\C^\mathrm{T}}.
$$
On the other hand, note that
the counit of $C^{\cop}$ is also $\varepsilon$, which satisfies
$\varepsilon(\C^\mathrm{T})= I_r$. Consequently, $\C^\mathrm{T}$ is a multiplicative matrix over $C^\cop$.
%
%By $\Delta(c_{ij}\otimes d_{st})=\sum_{k=1}^r
%\sum_{l=1}^s (c_{ik}\otimes d_{sl})\otimes (c_{kj}\otimes d_{lt})$ and definition of multiplication matrix, we obtain $\C\;\widetilde{\otimes}\;\D$ is the multiplication matrix of $C\otimes D$.
\end{proof}

Now we introduce the following notations: Let $U$ and $V$ be vector spaces.
For matrices $\A=(u_{ij})_{m\times n}$ over $U$ and $\B=(v_{ij})_{l\times p}$ over $V$, we denote
$$
u\otimes\B=
\left(\begin{array}{cccc}
    u\otimes v_{11} & u\otimes v_{12} & \cdots & u\otimes v_{1s}  \\
    u\otimes v_{21} & u\otimes v_{22} & \cdots & u\otimes v_{2s}  \\
    \vdots & \vdots & \ddots & \vdots  \\
    u\otimes v_{s1} & u\otimes v_{s2} & \cdots & u\otimes v_{ss}
\end{array}\right)\;\;\;\;(\forall u\in U),
$$
and define the \textit{Kronecker tensor product}
$\A\otimes\B$ to be the $ml\times np$ matrix over $U\otimes V$ given by
$$
\A\otimes\B=
\left(\begin{array}{cccc}
    u_{11}\otimes\B & u_{12}\otimes\B & \cdots & u_{1r}\otimes\B  \\
    u_{21}\otimes\B & u_{22}\otimes\B & \cdots & u_{2r}\otimes\B  \\
    \vdots & \vdots & \ddots & \vdots  \\
    u_{r1}\otimes\B & u_{r2}\otimes\B & \cdots & u_{rr}\otimes\B
\end{array}\right).
$$

From the above definition of the Kronecker tensor product, the following lemma is  straightforward to verify.

\begin{lemma}\label{uIsB=uB}
Let $U$ and $V$ be vector spaces.
Then for any $u\in U$ and any matrix $\B=(v_{ij})_{l\times p}$ over $V$,
\begin{equation*}
(uI_s)\;\widetilde{\otimes}\;\B= u\otimes\B\in\mathbf{M}_{l\times p}(U\otimes V).
\end{equation*}
\end{lemma}

\begin{proposition}\label{otimes}
Let $\C=(c_{ij})_{r\times r}$ and $\D=(d_{ij})_{s\times s}$ be multiplicative matrices over coalgebras $H$ and $K$, respectively. Then the following statements hold.
\begin{itemize}
\item[(1)]
The matrix $\C\otimes\D$ is multiplicative over $H\otimes K$.
\item[(2)]
For any positive integer $t$,
the matrices
$$
\C\odot I_t=
\left(\begin{array}{cccc}
    c_{11}I_t & c_{12}I_t & \cdots & c_{1r}I_t  \\
    c_{21}I_t & c_{22}I_t & \cdots & c_{2r}I_t  \\
    \vdots & \vdots & \ddots & \vdots  \\
    c_{r1}I_t & c_{r2}I_t & \cdots & c_{rr}I_t
\end{array}\right)\;\;\;\;
\text{and} \;\;\;\;
I_t\odot \D=
\left(\begin{array}{cccc}

    \D & 0 & \cdots & 0  \\
    0 & \D & \cdots & 0  \\
    \vdots & \vdots & \ddots & \vdots  \\
    0 & 0 & \cdots & \D
\end{array}\right)
$$
are multiplicative over $H$ and $K$, respectively.
\end{itemize}
\end{proposition}

\begin{proof}
\begin{itemize}
\item[(1)]
Since $\C$ and $\D$ are multiplicative,
we obtain that
$$
\Delta_{C\otimes D}(c_{ij}\otimes d_{kl})=\sum_{p=1}^r\sum_{q=1}^s
(c_{ip}\otimes d_{kq})\otimes(c_{pj}\otimes d_{ql})
$$
and
$$
\varepsilon_{C\otimes D}(c_{ij}\otimes d_{kl})=
\varepsilon_C(c_{ij})\varepsilon_D(d_{kl})=\delta_{ij}\delta_{kl}
$$
hold for all $1\leq i, j\leq r$ and $1\leq k, l\leq s$.
These equations can be concluded as the following equations
$$
\Delta_{C\otimes D}(\C\otimes\D)=(\C\otimes\D)\;\widetilde{\otimes}\;(\C\otimes\D)
\;\;\;\;\text{and}\;\;\;\;\varepsilon_{C\otimes D}(\C\otimes\D)=I_{rs}.
$$
In other words $\C\otimes\D$ is multiplicative over $H\otimes K$.
\item[(2)]
Clearly, the identity matrix is multiplicative over $\k$, since $\Delta(1)=1\otimes1$ and $\varepsilon(1)=1$.
%Applying (1) to the pairs $(\C, I_t)$ and $(I_{t}, \D)$, and using the canonical coalgebra isomorphisms
%    $$
%    H\otimes\k\cong H,\;\;\;\; \k\otimes K\cong K,
%    $$
%    we obtain the multiplicative matrices $\C\odot I_t$ over $H$ and $I_t\odot\D$ over $K$, as required.
Then we know by (1) that $\C\otimes I_t$ is a multiplicative matrix over $H\otimes K$. However, its image under the canonical coalgebra isomorphism $H\otimes\k\cong H$ will become $\C\odot I_t$, which is hence also multiplicative. A similar argument shows that $I_t\odot D$ is multiplicative over $K$.
\end{itemize}
\end{proof}

Finally, these constructions also imply the following results on basic multiplicative matrices.

\begin{corollary}\label{CT}
Suppose $\C$ and $\D$ are basic multiplicative matrices over $C$ and $D$, respectively. Then
\begin{itemize}
\item[(1)]
$\C^\mathrm{T}$ is a basic multiplicative matrix over $C^{\cop}$.
\item[(2)]
$\C\otimes\D$ is a basic multiplicative matrix over $C\otimes D$.
\end{itemize}
\end{corollary}

\section{Matrix basis decomposition of graded coalgebras}

In this section, we still use the notation in Section 2. For a matrix $\mathcal{A}\in\mathbf{M}_{r\times s}(U)$, we denote by $\mathrm{span}(\mathcal{A})$ the subspace spanned by its entries of $U$.

\subsection{Matrix descriptions for bases of simple (bi)comodules}

We begin with a matrix description of simple comodules.

\begin{lemma}\label{rhoM=CM}
Let $C$ be a simple coalgebra.
For any  basic multiplicative matrix $\C=(c_{ij})_{r\times r}$ of $C$ and any simple left $C$-comodule $(M,\rho)$, there exists a matrix $\M\in\mathbf{M}_{r\times1}(M)$ such that
$$\rho(\mathcal{M})=\C\;\widetilde{\otimes}\;\mathcal{M}$$
and the entries of $\mathcal{M}$ form a basis of $M$.
\end{lemma}

\begin{proof}
Since $\C=(c_{ij})_{r\times r}$ is a multiplicative matrix, one can find
$$
\Delta
\left(\begin{array}{c}
    c_{11} \\
    c_{21} \\
    \vdots  \\
    c_{r1}
\end{array}\right)
=\C\;\widetilde{\otimes}\;
\left(\begin{array}{c}
    c_{11} \\
    c_{21} \\
    \vdots  \\
    c_{r1}
\end{array}\right).
$$
Thus $V:=\k\{c_{i1}\mid1\leq i\leq r\}$ is a left $C$-subcomodule of $C$, and $\dim(V)=r$ as $\C$ is basic.
It follows that $(V, \Delta)$ is a simple left $C$-comodule.
Then under the isomorphism $(V, \Delta)\cong(M, \rho)$ of $C$-comodules, the image $\M$ of $(c_{11}, c_{21}, \cdots, c_{r1})^\mathrm{T}$ satisfies $$\rho(\mathcal{M})=\C\;\widetilde{\otimes}\;\mathcal{M}.$$
\end{proof}

The case of bicomodule follows by combining Lemmas \ref{cong} and \ref{rhoM=CM}.
%to the tensor product coalgebra $C\otimes D^{\cop}$.

\begin{proposition}\label{MCDBMM}
Let $C$ and $D$ be simple coalgebras with basic multiplicative matrices $\mathcal{C}_{r\times r}$ and $\mathcal{D}_{s\times s}$, respectively. For any simple $C$-$D$-bicomodule $(M, \rho^L, \rho^R)$, there exists a matrix $\M\in\mathbf{M}_{r\times s}(M)$
such that $$\rho^L(\mathcal{M})=\mathcal{C}\;\widetilde{\otimes}\;\mathcal{M}
,\;\;\;\;
\rho^R(\mathcal{M})=\mathcal{M}\;\widetilde{\otimes}\;\mathcal{D},
$$
and the entries of $\M$ form a basis of $M$.
\end{proposition}

\begin{proof}
By Lemma \ref{cong}, the isomorphism $\Phi$ sends $(M, \rho^L, \rho^R)$ to a simple left $C\otimes D^{\cop}$-comodule $(M, \delta)$. Meanwhile, we know by Corollary \ref{CT} implies that
%and Remark \ref{CT}
$\C\otimes\D^\mathrm{T}$ is a basic multiplicative matrix of $C\otimes D^{\cop}$.
Thus according to
Lemma \ref{rhoM=CM},  there exists a basis $\{m_{ij}\mid1\leq i\leq r, 1\leq j\leq s\}$ of $M$, such that the matrix
$$\mathcal{M^\wr}=
\left(\begin{array}{r}
m_{11}, m_{12}, \cdots, m_{1s},\;  m_{21}, m_{22}, \cdots, m_{2s},\;
\cdots,\; m_{r1}, m_{r2}, \cdots, m_{rs}
\end{array}\right)^\mathrm{T}
$$
satisfies
$\delta({\M}^\wr)=(\C\otimes \D^\mathrm{T})\;\widetilde{\otimes}\;{\M}^\wr.$
%Equivalenty, $\delta(m_{ij})=\sum_{k=1}^r\sum_{l=1}^s(c_{ik}\otimes d_{lj})\otimes m_{kl}$.

Let us set $\M=(m_{ij})_{r\times s}$ and write $\M_i=(m_{i1}, m_{i2}, \cdots, m_{is})^\mathrm{T}$ for $1\leq i\leq r$.
Thus
$$
\mathcal{M^\wr}=
\left(\begin{array}{c}
\M_1 \\ \M_2 \\ \vdots \\ \M_r
\end{array}\right)\;\;\;\;\text{and}\;\;\;\;
\mathcal{M}=
\left(\begin{array}{c}
\M_1^\mathrm{T} \\ \M_2^\mathrm{T} \\ \vdots \\  \M_r^\mathrm{T}
\end{array}\right).
$$

Now consider $(M, \rho^L, \rho^R)$ to be the preimage of $(M, \delta)$ under $\Phi$. Recall from the proof of Lemma \ref{cong} that $\rho^L$ and $\rho^R$ should satisfy Equations (\ref{rhoL}) and (\ref{rhoR}). Consequently, since $\varepsilon_D(\D^\mathrm{T})=I_s$, we have
\begin{eqnarray*}
\left(\begin{array}{c}
\rho^L(\M_1) \\ \rho^L(\M_2) \\ \vdots \\ \rho^L(\M_r)
\end{array}\right)
&=&
\rho^L(\M^\wr)=[\C\odot\varepsilon(\D^\mathrm{T})]
  \;\widetilde{\otimes}\;\mathcal{M^\wr} \\
&=&
\left(\begin{array}{cccc}
    c_{11}I_s & c_{12}I_s & \cdots & c_{1r}I_s  \\
    c_{21}I_s & c_{22}I_s & \cdots & c_{2r}I_s \\
    \vdots & \vdots & \ddots & \vdots  \\
    c_{r1}I_s & c_{r2}I_s & \cdots & c_{rr}I_s
\end{array}\right)
\;\widetilde{\otimes}\;
\left(\begin{array}{c}
\M_1 \\ \M_2 \\ \vdots \\ \M_r
\end{array}\right)
=
\left(\begin{array}{c}
    \sum_{k=1}^r c_{1k}\otimes \M_k \\
    \sum_{k=1}^r c_{2k}\otimes \M_k \\
    \vdots   \\
    \sum_{k=1}^r c_{rk}\otimes \M_k
\end{array}\right),
\end{eqnarray*}
where the last equality is due to Lemma \ref{uIsB=uB}.
It follows that $\rho^L(\M_i)=\sum_{k=1}^rc_{ik}\otimes \M_k$, and hence $\rho^L(\M_i)^\mathrm{T}=\sum_{k=1}^r c_{ik}\otimes\M_k^\mathrm{T}$ for each $i$. Then
\begin{eqnarray*}
\left(\begin{array}{c}
\rho^L(\M_1)^\mathrm{T} \\ \rho^L(\M_2)^\mathrm{T} \\
\vdots \\ \rho^L(\M_r)^\mathrm{T}
\end{array}\right)
=
\left(\begin{array}{c}
\sum_{k=1}^r c_{1k}\otimes \M_k^\mathrm{T} \\
\sum_{k=1}^r c_{2k}\otimes \M_k^\mathrm{T} \\
\vdots \\
\sum_{k=1}^r c_{rk}\otimes \M_k^\mathrm{T} \\
\end{array}\right)
=
\left(\begin{array}{cccc}
c_{11} & c_{12} & \cdots & c_{1r} \\
c_{21} & c_{22} & \cdots & c_{2r} \\
\vdots & \vdots & \ddots & \vdots \\
c_{r1} & c_{r2} & \cdots & c_{rr}
\end{array}\right)
\;\widetilde{\otimes}\;
\left(\begin{array}{c}
\M_1^\mathrm{T} \\ \M_2^\mathrm{T} \\
\vdots \\ \M_r^\mathrm{T}
\end{array}\right).
\end{eqnarray*}
In other words, $\rho^L(\M)=\C\;\widetilde{\otimes}\;\M$.

Similarly, since $\varepsilon_C(\C)=I_r$, the right coaction $\rho^R$ satisfies
$$\rho^R(\M^\wr)=\tau_{D, M}[(I_r\otimes\D^\mathrm{T})\;\widetilde{\otimes}\;
\mathcal{M^\wr}]=\tau_{D, M}[(I_r\odot\D^\mathrm{T})\;\widetilde{\otimes}\;\M^\wr].
$$
However, we have
\begin{eqnarray*}
(I_r\odot\D^\mathrm{T})\;\widetilde{\otimes}\;
\mathcal{M^\wr}
&=&
\left(\begin{array}{cccc}
    \D^\mathrm{T} & 0 & \cdots & 0 \\
    0 &  \D^\mathrm{T} & \cdots & 0 \\
    \vdots& \vdots & \ddots & \vdots \\
    0& 0 &\cdots & \D^\mathrm{T}
\end{array}\right)
\;\widetilde{\otimes}\;
\left(\begin{array}{c}
\M_1\\  \M_2\\ \vdots \\ \M_r
\end{array}\right)
=
\left(\begin{array}{c}
    \D^\mathrm{T}\;\widetilde{\otimes}\;\M_1  \\
    \D^\mathrm{T}\;\widetilde{\otimes}\;\M_2 \\
    \vdots   \\
    \D^\mathrm{T}\;\widetilde{\otimes}\;\M_r
\end{array}\right). \\
%&=&
%\left(\begin{array}{c}
%    m_{11}\otimes d_{11}+ m_{12}\otimes d_{21}+\cdots+  m_{1s}\otimes d_{s1}  \\
%    \vdots \\
%    m_{11}\otimes d_{1s}+ m_{12}\otimes d_{2s}+\cdots+  m_{1s}\otimes d_{ss} \\
%    \vdots   \\
%    m_{r1}\otimes d_{11}+ m_{r2}\otimes d_{21}+\cdots+  m_{rs}\otimes d_{s1}  \\
%    \vdots \\
%    m_{r1}\otimes d_{s1}+ m_{r2}\otimes d_{2s}+\cdots+  m_{rs}\otimes d_{ss}
%\end{array}\right)
%=
%\left(\begin{array}{c}
%\rho^R(m_{11}) \\ \rho^R(m_{12})\\ \vdots\\\rho^R(m_{r(s-1)})\\ \rho^R(m_{rs})
%\end{array}\right)
\end{eqnarray*}
Also, one can verify that
$$
\tau_{D, M}(\D^\mathrm{T}\;\widetilde{\otimes}\;\M_i)=
(\M_i^\mathrm{T}\;\widetilde{\otimes}\;\D)^\mathrm{T}
$$
holds for each $i$.

Then by combining these equations, we obtain
\begin{eqnarray*}
\left(\begin{array}{c}
\rho^R(\M_1) \\ \rho^R(\M_2) \\ \vdots \\
\rho^R(\M_r)
\end{array}\right)
=
\rho^R(\M')=
\left(\begin{array}{c}
\tau_{D, M}(\D^\mathrm{T}\;\widetilde{\otimes}\;\M_1) \\
\tau_{D, M}(\D^\mathrm{T}\;\widetilde{\otimes}\;\M_2) \\
\vdots \\
\tau_{D, M}(\D^\mathrm{T}\;\widetilde{\otimes}\;\M_r)
\end{array}\right)
=
\left(\begin{array}{c}
(\M_1^\mathrm{T}\;\widetilde{\otimes}\;\D)^\mathrm{T} \\
(\M_2^\mathrm{T}\;\widetilde{\otimes}\;\D)^\mathrm{T} \\
\vdots \\
(\M_r^\mathrm{T}\;\widetilde{\otimes}\;\D)^\mathrm{T}
\end{array}\right).
\end{eqnarray*}
%Consequently, $\rho^R(\M_i)^\mathrm{T}= (\rho^R(\M_i))^\mathrm{T} = \M_i^\mathrm{T}\;\widetilde{\otimes}\;\D$ for $1\leq i\leq r$.
Finally, we conclude that
\begin{eqnarray*}
\rho^R(\M)=
\left(\begin{array}{c}
\rho^R(\M_1)^\mathrm{T}\\
\rho^R(\M_2)^\mathrm{T} \\
\vdots \\
\rho^R(\M_r)^\mathrm{T}
\end{array}\right)
=
\left(\begin{array}{c}
\M_1^\mathrm{T}\;\widetilde{\otimes}\;\D \\
\M_2^\mathrm{T}\;\widetilde{\otimes}\;\D \\
\vdots \\
\M_r^\mathrm{T}\;\widetilde{\otimes}\;\D
\end{array}\right)
=
\left(\begin{array}{c}
\M_1^\mathrm{T} \\
\M_2^\mathrm{T} \\
\vdots \\
\M_r^\mathrm{T}
\end{array}\right)\;\widetilde{\otimes}\;\D
=\M\;\widetilde{\otimes}\;\D.
\end{eqnarray*}
\end{proof}

\begin{lemma}\label{N=M}
Let ${C}$ be a simple coalgebra with a basic multiplicative matrix $\C_{r\times r}$, and let $(V, \delta)\in {}^C\mathfrak{M}$.
Suppose that $V=\bigoplus_{\mathcal{M}\in\mathbf{\Gamma}}\mathrm{span}
(\mathcal{M})$, where
$\mathbf{\Gamma}\subseteq\mathbf{M}_{r\times 1}(V)$ is a family of nonzero matrices satisfying
$$\delta(\M)=\C\;\widetilde{\otimes}\;\M\;\;\;\;(\M\in\mathbf{\Gamma}).$$
If $\mathcal{N}\in\mathbf{M}_{r\times1}(V)$ satisfies $\delta(\mathcal{N})\in\mathbf{M}_{r\times1}(V)$ and $\delta(\N)=\C\;\widetilde{\otimes}\;\N$, then there exist unique scalars $\alpha_{\mathcal{M}}\in\k\;(\M\in\mathbf{\Gamma})$, all but finitely many zero, such that $ \mathcal{N}=\sum_{\mathcal{M}\in\mathbf{\Gamma}}
\alpha_\mathcal{M}\mathcal{M}$.
\end{lemma}

\begin{proof}
Firstly, for each $\M\in\mathbf{\Gamma}$, it is evident that $\mathrm{span}(\M)\in{}^C\mathfrak{M}$ with $0<\dim(\mathrm{span}(\M))\leq r$, which implies that $\mathrm{span}(\M)$ is a simple left $C$-comodule of dimension $r$. As a consequence, the entries of all $\M\in\mathbf{\Gamma}$ form a basis of $V=\bigoplus_{\mathcal{M}\in\mathbf{\Gamma}}\mathrm{span}
(\mathcal{M})$.

Now we write $\mathcal{N}=(v_1, v_2, \cdots, v_r)^\mathrm{T} \in\mathbf{M}_{r\times 1}(V)$.
Then
for each $1\leq i\leq r$ there exist unique row vectors $A_i^{(\M)}\in \mathbf{M}_{1\times r}(\k)\;(\M\in\mathbf{\Gamma})$, only finitely many nonzero, such that $v_i=\sum_{\M\in\mathbf{\Gamma}}A_i{}^{(\M)}\M.$
Thus, we obtain
$$\mathcal{N}_{r\times1}=
\left(\begin{array}{c}
    \sum_{\mathcal{M}\in\mathbf{\Gamma}}A_1^{(\mathcal{M})}
    \mathcal{M} \\
    \sum_{\mathcal{M}\in\mathbf{\Gamma}}A_2^{(\mathcal{M})}
    \mathcal{M}  \\
    \vdots  \\
    \sum_{\mathcal{M}\in\mathbf{\Gamma}}A_r^{(\mathcal{M})}
    \mathcal{M}
\end{array}\right)
=\sum_{\mathcal{M}\in\mathbf{\Gamma}}
\left(\begin{array}{c}
A_1^{(\mathcal{M})}\\
A_2^{(\mathcal{M})}  \\
\vdots  \\
A_r^{(\mathcal{M})}
\end{array}\right)
\mathcal{M}.
$$
For each $\M\in\mathbf{\Gamma}$, set the $r\times r$ matrix
$$
A^{(\M)}=\left(\begin{array}{c}
A_1^{(\M)}\\
A_2^{(\M)}\\
\vdots \\
A_r^{(\M)}
\end{array}\right)\in\mathbf{M}_r(\k).
$$
Then $\mathcal{N}=\sum_{\M\in\mathbf{\Gamma}}A^{(\M)}\M$, where only finitely many terms are nonzero.
It follows that
$$
\delta(\mathcal{N})=\C\;\widetilde{\otimes}\;\mathcal{N}=
\sum_{\mathcal{M}\in\mathbf{\Gamma}}\C\;\widetilde{\otimes}\;
(A^{(\mathcal{M})}\mathcal{M})=
\sum_{\mathcal{M}\in\mathbf{\Gamma}}(\C A^{(\mathcal{M})})
\;\widetilde{\otimes}\;\mathcal{M}.
$$
Meanwhile, we can find $\delta({\mathcal{N}})=\sum_
{\mathcal{M}\in\mathbf{\Gamma}}A^{(\mathcal{M})}\delta(\mathcal{M})=
\sum_{\mathcal{M}\in\mathbf{\Gamma}}(A^{(\mathcal{M})}\C)
\;\widetilde{\otimes}\;\mathcal{M}$.
Thus,
$$
\sum_{\M\in\mathbf{\Gamma}}(\C A^{(\M)}-A^{(\M)}\C)\;\widetilde{\otimes}\;\M=0.
$$
However, since the entries of all $\M\in\mathbf{\Gamma}$ are linearly independent,
we know by Remark \ref{AB}(2) that holds $\C A^{(\mathcal{M})}=A^{(\mathcal{M})}\C$ for each $\mathcal{M}\in\mathbf{\Gamma}$.
Consequently, one can find that $A^{(\mathcal{M})}=\alpha_{\mathcal{M}}I_r$ holds for some $\alpha_{\mathcal{M}}\in\k$, because the entries of $\C$ are linearly independent.
\end{proof}

Similarly, according to the isomorphism $\Phi: {}^C\mathfrak{M}{}^D\cong {}^{C\otimes D^\cop}\mathfrak{M}$ in Lemma \ref{cong}, one can obtain the bicomodule version of Lemma \ref{N=M} as the following corollary.

\begin{corollary}\label{NM}
Let $C$ and $D$ be simple coalgebras with basic multiplicative matrices $\C_{r\times r}$ and $\D_{s\times s}$, respectively, and let $(V, \rho^L, \rho^R)\in {}^C\mathfrak{M}{}^D$. Suppose that
$V=\bigoplus_{\mathcal{M}\in\mathbf{\Gamma}}
\mathrm{span}(\mathcal{M})$, where
$\mathbf{\Gamma}\subseteq\mathbf{M}_{r\times s}(V)$ is a family of nonzero matrices satisfying  $$\rho^L(\mathcal{M})=\C\;\widetilde{\otimes}\;
\mathcal{M},\;\;\;\;
\rho^R(\mathcal{M})=\mathcal{M}
\;\widetilde{\otimes}\;\D\;\;\;\;(\mathcal{M}\in\mathbf{\Gamma}).$$
If $\mathcal{N}\in\mathbf{M}_{r\times s}(V)$ satisfies $$\rho^L(\mathcal{N})=\C
\;\widetilde{\otimes}\;\mathcal{N},\;\;\;\;
\rho^R(\mathcal{N})=\mathcal{N}\;\widetilde{\otimes}\;\D,$$  then there exist unique scalars $\alpha_{\mathcal{M}}\in\k\;(\M\in\mathbf{\Gamma})$,
all but finitely many zero, such that $ \mathcal{N}=\sum_{\mathcal{M}\in\mathbf{\Gamma}}
\alpha_\mathcal{M}\mathcal{M}$.
\end{corollary}

\begin{proof}
%Let $\delta: V\rightarrow(C\otimes D^{\cop})\otimes V$ be the left coaction corresponding to the $C$-$D$-bicomodule structure on $V$ by Lemma \ref{cong}. For each $\M=(m_{ij})\in\mathbf{\Gamma}$, set
%$$
%\M'=
%\left(\begin{array}{r}
%m_{11}, m_{12}, \cdots  m_{1s}, \; m_{21}, \cdots m_{2s},
%\cdots m_{r1}, \; m_{r2}, \cdots m_{rs}
%\end{array}\right)^\mathrm{T},
%$$
%and define $\N'$ analogously.
%Then, we have
Let us use the following notation:
For any vector space $U$, and any matrix $\A\in \mathbf{M}_{r\times s}(U)$, we denote its vectorization by
$$
\A^\wr=
\left(\begin{array}{r}
u_{11}, u_{12}, \cdots, u_{1s}, \; u_{21}, u_{22}, \cdots, u_{2s},\;
\cdots,\; u_{r1}, u_{r2}, \cdots, u_{rs}
\end{array}\right)^\mathrm{T}\in\mathbf{M}_{rs\times1}(U),
$$

Now consider $(V, \delta)=\Phi(V, \rho^L, \rho^R)$. A similar argument with the proof of Proposition \ref{MCDBMM} provides that
$$
\delta(\M^\wr)=(\C\otimes\D^\mathrm{T})
\;\widetilde{\otimes}\;\M^\wr\;\;\;\;(\forall \M\in\mathbf{\Gamma}),
\;\;\;\;
\delta(\N^\wr)=(\C\otimes\D^\mathrm{T})\;\widetilde{\otimes}\;\N^\wr.
$$
By Corollary \ref{CT}, $\C\otimes\D^\mathrm{T}$ is a basic multiplicative matrix of the simple coalgebra $C\otimes D^{\cop}$.
%Moreover, $\M'\neq0$ for every $\M\in\mathbf{\Gamma}$, and $V=\bigoplus\limits_{\mathcal{M}\in\mathbf{\Gamma}}
%\mathrm{span}(\mathcal{M}')$.
Then according to Lemma \ref{N=M}, we obtain unique scalars $\alpha_{\M}\in\k$, all but finitely many zero, such that
$
\mathcal{N}^\wr=\sum_{\mathcal{M}\in\mathbf{\Gamma}}
\alpha_\mathcal{M}\mathcal{M}^\wr.
$
Namely,
$
\mathcal{N}=\sum_{\mathcal{M}\in\mathbf{\Gamma}}
\alpha_\mathcal{M}\mathcal{M}
$.
\end{proof}

The next lemma describes a basis for the cotensor product of two simple bicomudules.

\begin{lemma}\label{MN=span}
Let $C$, $D$ and $E$ be simple coalgebras with basic multiplicative matrices $\C_{r\times r}$, $\D_{s\times s}$ and $\E_{t\times t}$, respectively. Let $(M, \rho^L_M, \rho^R_M)$ be a simple $C$-$E$-bicomodule and $(N, \rho^L_N, \rho^R_N)$ a simple $E$-$D$-bicomodule.
Suppose that the entries of $\M\in\mathbf{M}_{r\times t}(M)$ and $\N\in\mathbf{M}_{t\times s}(N)$ form bases of $M$ and $N$, respectively, and that
$$
\rho_{M}^L(\M)=\C\;\widetilde{\otimes}\;\M,\;\;\;\;
\rho_{M}^R(\M)=\M\;\widetilde{\otimes}\;\E,\;\;\;\;
\rho_{N}^L(\N)=\E\;\widetilde{\otimes}\;\N,\;\;\;\;
\rho_{N}^R(\N)=\N\;\widetilde{\otimes}\;\D.
$$
Then the entries of $\M\;\widetilde{\otimes}\;\N$ form a basis of $M\square_E N$.

%If the cotensor product $M\square_EN$ is simple as a $C$-$D$-bicomodule, then
%\begin{itemize}
%\item[\textup{(1)}]
%$M$ and $N$ are simple as a $C$-$E$-bicomodule and an $E$-$D$-bicomodule, respectively.
%\item[\textup{(2)}] there exist matrices
%$\M$, $\mathcal{N}$,
%whose entries form bases of $M$ and $N$, respectively, such that
%$$
%\rho_{M}^L(\M)=\C\;\widetilde{\otimes}\;\M,\;\;
%\rho_{M}^R(\M)=\M\;\widetilde{\otimes}\;\E,\;\;
%\rho_{N}^L(\N)=\E\;\widetilde{\otimes}\;\N,\;\;
%\rho_{N}^R(\N)=\N\;\widetilde{\otimes}\;\D.
%$$
%Moreover, the entries of $\M\;\widetilde{\otimes}\;\N$
%form a basis of $M\square_EN$. In particular, $M\square_EN=span_{\k}(\M\;\widetilde{\otimes}\;\N)$.
%\end{itemize}
\end{lemma}

\begin{proof}
According to Corollary \ref{dim=rs}, it is clear that
$\dim(M)=rt$ and $\dim(N)=ts$. Then by Lemma \ref{cotenp dim} to the right $E$-comodule and the left $E$-comodule $N$, we obtain
\begin{eqnarray*}
\dim(M\square_E N)=\frac{\dim(M)\cdot\dim(N)}{\dim(E)}=
\frac{rt\cdot ts}{t^2}=rs.
\end{eqnarray*}
%Next, we consider the tensor product basis of $M\otimes N$ induced by the entries of $\M$ and $\N$. By the definition of $\widetilde{\otimes}$, each entry of $\M\;\widetilde{\otimes}\;\N$ is a nonempty sum of distinct elements of this basis, and different entries involve disjoint sets of basis elements. Hence the entries of $\M\;\widetilde{\otimes}\;\N$ are linearly independent.
%

Moreover, it follows from the assumptions that
$$
(\rho^R_M\otimes\id_N)(\M\;\widetilde{\otimes}\;\N)=
(\M\;\widetilde{\otimes}\;\E)\;\widetilde{\otimes}\;\N
=\M\;\widetilde{\otimes}\;(\E\;\widetilde{\otimes}\;\N)
=(\id_M\otimes\rho^L_N)(\M\;\widetilde{\otimes}\;\N)
$$
hold, which implies that
$\M\;\widetilde{\otimes}\;\N\in\mathbf{M}_{r\times s}(M\square_EN)$.
Thus it suffices to explain that $\M\;\widetilde{\otimes}\;\N$ has linearly independent entries, and this fact is mentioned in Remark \ref{AB}(1).
\end{proof}

\subsection{A matrix decomposition for $\mathbb{N}$-graded coalgebra}
In this subsection, we always assume that $H=\bigoplus_{n\geq0}H(n)$ is an $\mathbb{N}$-graded coalgebra with $H(0)=H_0$. Let $\pi: H\rightarrow H(0)$ be the canonical projection of coalgebras. Then each $H(n)$ is an $H_0$-bicomodule via
\begin{equation}\label{rhoLR}
\rho^L=(\pi\otimes\id)\circ\Delta|_{H(n)}: H(n)\rightarrow H_0\otimes H(n),\;\;\;\;
\rho^R=(\id\otimes \pi)\circ\Delta|_{H(n)}:
H(n)\rightarrow H(n)\otimes H_0.
\end{equation}

Recall that $\mathcal{S}$ denotes the set of all the simple subcoalgebras of $H$, and
we use the notation ${}^CH(n){}^D$ for $C, D\in\mathcal{S}$ from Proposition \ref{Hn}.

Note that ${}^CH(n){}^D$ is a $C$-$D$-bicomodule with structures $\rho^L$ and $\rho^R$ (\ref{rhoLR}). In fact, one can find by Proposition \ref{Hn}(1) and (2) that $\rho^L({}^CH(n){}^D)\subseteq C\otimes {}^CH(n){}^D$ and $\rho^R({}^CH(n){}^D)\subseteq {}^CH(n){}^D\otimes D$ hold.

\begin{lemma}\label{DWcotp}
Let $C, D\in\mathcal{S}$ with basic multiplicative matrices $\C_{r\times r}$ and $\D_{s\times s}$, respectively. Let $n\geq1$ and $W\subseteq{}^CH(n){}^D$ is a simple $C$-$D$-subbicomodule. If $\W\in\mathbf{M}_{r\times s}(W)$ satisfies that
\begin{eqnarray}\label{rhoLRW}
\rho^L(\W)=\C\;\widetilde{\otimes}\;\W,\;\;\;\; \rho^R(\W)=\W\;
\widetilde{\otimes}\;\D,
\end{eqnarray}
and the entries of $\W$ form a basis of $W$. Then
\begin{eqnarray*}
\Delta(\W)-\C\;\widetilde{\otimes}\;\W
-\W\;\widetilde{\otimes}\;\D\in\bigoplus\limits_{i=1}^{n-1}
\bigoplus\limits_{E\in\mathcal{S}}
\mathbf{M}_{r\times s}
({}^CH(i){}^E\square_E {}^EH(n-i)^D).
\end{eqnarray*}
\end{lemma}

\begin{proof}
%Recall that $\rho^L=(\pi\otimes\id)\circ\Delta$ and $\rho^R=(\id\otimes \pi)\circ\Delta$.

Since the entries of $\W$ belong to ${}^CH(n){}^D$, it follows from Proposition \ref{Hn} (1) and (2) that
\begin{eqnarray*}\label{delta(W)in}
&&  \Delta(\W)
~\in~ \bigoplus\limits_{i=0}^n\bigoplus\limits_{E\in\mathcal{S}}
\mathbf{M}_{r\times s}({}^CH(i)^E\otimes {}^EH(n-i){}^D) \nonumber \\
&=&
\mathbf{M}_{r\times s}(C\otimes {}^CH(n){}^D)
\oplus\mathbf{M}_{r\times s}({}^CH(n){}^D\otimes D)
\oplus
\left[\bigoplus\limits_{i=1}^{n-1}\bigoplus\limits_{E\in\mathcal{S}}
\mathbf{M}_{r\times s}({}^CH(i)^E\otimes {}^EH(n-i){}^D)\right].
%&=&
%\big[\C\;\widetilde{\otimes}\;\W+
%\W\;\widetilde{\otimes}\;\D+
%\bigoplus\limits_{i=0}^{n-1}\bigoplus\limits_{E\in\mathcal{S}}
%\mathbf{M}_{r\times s}({}^CH(i)^E\otimes {}^EH(n-i){}^D)\big].
\end{eqnarray*}
Then applying $\pi\otimes\id$ and $\id\otimes\pi$ respectively on the both sides, we obtain by (\ref{rhoLRW}) that
\begin{eqnarray*}
\Delta(\W)-\C\;\widetilde{\otimes}\;\W-\W\;\widetilde{\otimes}\;\D
\in\bigoplus\limits_{i=1}^{n-1}\bigoplus_{E\in\mathcal{S}}
\mathbf{M}_{r\times s}({}^CH(i){}^E\otimes{}^E H(n-i)^D).
\end{eqnarray*}
%In fact, we show that $\Delta(\W)-\C\;\widetilde{\otimes}\;\W
%-\W\;\widetilde{\otimes}\;\D\in\bigoplus\limits_{i=0}^n
%\mathcal{M}_{r\times s}
%({}^CH(i)\square_{H_0} H(n-i)^D$,
Furthermore, it is straightforward to find that the images of $\Delta(\W)-\C\;\widetilde{\otimes}\;
\W-\W\;\widetilde{\otimes}\;\D$ under $\rho^R\otimes\id$ and $\id\otimes\rho^L$ are equal. In fact, the image is $(\id\otimes\pi\otimes\id)\circ\Delta_3(\W)-\C\;\widetilde{\otimes}\;
\C\;\widetilde{\otimes}\;\W-\W\;\widetilde{\otimes}\;
\D\;\widetilde{\otimes}\;\D$, where $\Delta_3=(\Delta\otimes\id)\circ\Delta.$
Consequently,
\begin{eqnarray*}
\Delta(\W)-\C\;\widetilde{\otimes}\;\W-\W\;\widetilde{\otimes}\;\D
\in\bigoplus\limits_{i=1}^{n-1}\bigoplus_{E\in\mathcal{S}}
\mathbf{M}_{r\times s}({}^CH(i){}^E\square_E{}^E H(n-i)^D)
\end{eqnarray*}
holds since the cotensor product $-\square_E-$ is biadditive.
\end{proof}

We now combine the preceding results to
describe matrix bases for all homogeneous components of $H$ and describe their coproducts.

\begin{theorem}\label{mainthem}
Fix a basic multiplicative matrix for each simple subcoalgebra of $H$, and let $\mathcal{P}_0$ be the set of all such matrices. Then for each $n\geq1$, there exist a set $\mathcal{P}_{n}$ of matrices over $H(n)$ satisfying the following properties:
\begin{itemize}
\item[(1)]  Every matrix in $\mathcal{P}_n$ has linearly independent entries.
\item [(2)] For any $\C, \D\in \mathcal{P}_0$, let $C=\mathrm{span}(\C)$ and $D=\mathrm{span}(\D)$. There exists a subset ${}^\C{\mathcal{P}_n}{}^\D\subseteq\mathcal{P}_n$ such that
    $$
    {}^C{H(n)}{}^D=\bigoplus_{\mathcal{W}
    \in{}^\C{\mathcal{P}_n}^\D}\mathrm{span}(\mathcal{W}).
    $$
\item [(3)] $\mathcal{P}_n$ is the disjoint union:
    $$
    \mathcal{P}_n=\bigcup\limits_{\C, \D\in\mathcal{P}_0}
    {}^\C{\mathcal{P}_n}{}^\D.
    $$
\item [(4)] For any $\C, \D\in \mathcal{P}_0$ and every $\mathcal{W}\in {}^\C{\mathcal{P}_n}^\D$, there exist $\alpha_{\mathcal{X},\mathcal{Y}}\in\k$, such that
    \begin{eqnarray}\label{deltaW}
    \Delta(\mathcal{W})=\C\;\widetilde{\otimes}\;\mathcal{W}
    +\mathcal{W}\;\widetilde{\otimes}\;\D
    +\sum_{i=1}^{n-1}\sum_{\mathcal{E}\in\mathcal{P}_0}
    \sum_{\substack{\mathcal{X}\in{}^\C{\mathcal{P}_i}^\mathcal{E},\\
    \mathcal{Y}\in{}^\mathcal{E}{\mathcal{P}_{n-i}}^\D}}
    \alpha_{\mathcal{X},\mathcal{Y}}(\mathcal{X}\;\widetilde{\otimes}\;
    \mathcal{Y}),
    \end{eqnarray}
where the coefficients $\alpha_{\mathcal{X},\mathcal{Y}}$ are uniquely determined.
\end{itemize}
\end{theorem}

\begin{proof}
We first construct the families $\mathcal{P}_n$ $(n\geq1)$ satisfying (1)-(3).
%Let $p_j: H\rightarrow H(j)$ denote the canonical projection.
%Since $H_0=H(0)$, and $\pi$ is a coalgebra map, each $H(n)$ has the $H_0$-bicomodule structure.
%We construct the families $\{\mathcal{P}_n\}_{n\geq1}$ inductively.
Since $C=\mathrm{span}(\C)$ and $D=\mathrm{span}(\D)$ are simple coalgebras, we know by Lemma \ref{cong} that the $C$-$D$-bicomodule ${}^CH(n){}^D$ is semisimple for $n\geq1$.
%and suppose that the families in all positive degrees less than $n$ have already been chosen to satisfy (1), (2), and the two coaction identities used in the construction.
%
%For $\C, \D\in\mathcal{P}_0$, put $C=\mathrm{span}(\C)$ and $D=\mathrm{span}(\D)$. We know that the $C$-$D$-bicomodule ${}^CH(n){}^D$ is  semisimple.
%Since $\C_{r\times r}$ and $\D_{s\times s}$ are basic multiplicative matrices,
%$$(C\otimes D^{\cop})^\ast\cong M_r(\k)\otimes M_s(\k)^{\op}\cong M_{rs}(\k).
%$$
%By Lemma \ref{cong},
%the category of ${}^C\mathfrak{M}{}^D$ is semisimple. Hence $H(n)$ admits a decomposition
Hence there is a decomposition
$${}^CH(n){}^D=\bigoplus_{\lambda\in\Lambda}W_{\lambda},$$ where each $W_{\lambda}$ is a simple $C$-$D$-subbicomodule. According to
Proposition $\ref{MCDBMM}$, for each $\lambda\in\Lambda$, there exists a matrix $\W_{\lambda}\in\mathbf{M}_{r\times s}(W_\lambda)$ such that
\begin{equation}\label{rhoLRWl}
\rho^L(\W_{\lambda})=\C\;\widetilde{\otimes}\;\W_{\lambda},\;\;
\;\;\rho^R(\W_{\lambda})=\W_{\lambda}\;\widetilde{\otimes}\;\D,
\end{equation}
and that the entries of $\W_\lambda$ form a basis of $W_\lambda$.

Set
$$
{}^\C\mathcal{P}_n{}^\D=\{\W_{\lambda}\mid \lambda\in\Lambda\},\;\;\;\;
\mathcal{P}_n=\bigcup\limits_{\C,\D\in
\mathcal{P}_0}{}^\C\mathcal{P}_n{}^\D.
$$
Consequently, one can easily check that all the properties in (1)-(3) hold. Note by (\ref{rhoLRWl}) that
\begin{equation}\label{rhoLRW1}
\rho^L(\W)=\C\;\widetilde{\otimes}\;\W\;\;\;\;\text{and}
\;\;\;\;\rho^R(\W)=\W\;\widetilde{\otimes}\;\D
\;\;\;\;(\forall\W\in{}^\C\mathcal{P}_n{}^\D).
\end{equation}

Next we aim to prove that (4) also holds for the above family $\{\mathcal{P}_n\}_{n\geq1}$.
Suppose $\C, \D, \E\in \mathcal{P}_0$,
and $1\leq i\leq n-1$,
%Set
%$$
%M=span_\k(\X),\;\;\;\;N=span_\k(\Y),
%$$
Then for any
$$\X\in{}^\C\mathcal{P}_i{}^\mathcal{E}, \;\;\;\; \Y\in{}^\mathcal{E}\mathcal{P}_{n-i}{}^\D,$$
we have by Lemma \ref{MN=span} that
$$
\mathrm{span}(\X)\square_E\mathrm{span}(\Y)=\mathrm{span}(\X\;\widetilde{\otimes}\;\Y),
$$
where $E=\mathrm{span}(\E)$. Furthermore, we obtain by the property in (2) for ${}^\C\mathcal{P}_i{}^\E$ and ${}^\E\mathcal{P}_{n-i}{}^\D$ that
\begin{eqnarray}\label{cotp=span}
\bigoplus\limits_{i=1}^{n-1}\bigoplus\limits_{E\in\mathcal{S}} {}^CH(i){}^E\square_E{}^EH(n-i)^D
&=&
\bigoplus\limits_{i=1}^{n-1}\bigoplus_{\mathcal{E}\in\mathcal{P}_0}
\big(\bigoplus\limits_{\X\in{}^\C{\mathcal{P}_i}{}^\mathcal{E}}
\mathrm{span}(\X)\big)\square_E
\big(\bigoplus\limits_{\Y\in{}^\mathcal{E}{\mathcal{P}_{n-i}}{}^\D}
\mathrm{span}(\Y)\big) \nonumber
\\
&=&
\bigoplus\limits_{i=1}^{n-1} \bigoplus_{\mathcal{E}\in\mathcal{P}_0}
\bigoplus\limits_{\substack{\X\in{}^\C{\mathcal{P}_i}^\mathcal{E},\\
    \mathcal{Y}\in{}^\mathcal{E}{\mathcal{P}_{n-i}}^\D}}
    (\mathrm{span}(\X)\square_E \mathrm{span}(\Y)) \nonumber \\
&=&
\bigoplus_{i=1}^{n-1} \bigoplus_{\mathcal{E}\in\mathcal{P}_0}
\bigoplus_{\substack{\X\in{}^\C{\mathcal{P}_i}^\mathcal{E},\\
    \mathcal{Y}\in{}^\mathcal{E}{\mathcal{P}_{n-i}}^\D}}   \mathrm{span}(\X\;\widetilde{\otimes}\;\Y),
\end{eqnarray}
where the second equality is because the cotensor product $-\square_E-$ is biadditive.

%$$
%(\rho^L\otimes\id)(\Delta(\W)-\C\;\widetilde{\otimes}\;
%\W-\W\;\widetilde{\otimes}\;\D)=\C\;\widetilde{\otimes}\;
%(\Delta(\W)-\C\;\widetilde{\otimes}\;
%\W-\W\;\widetilde{\otimes}\;\D),
%$$
%$$
%(\id\otimes\rho^R)(\Delta(\W)-\C\;\widetilde{\otimes}\;
%\W-\W\;\widetilde{\otimes}\;\D)=(\Delta(\W)-\C\;\widetilde{\otimes}\;
%\W-\W\;\widetilde{\otimes}\;\D)\;\widetilde{\otimes}\;\D.
%$$

Now suppose $\W\in{}^\C\mathcal{P}_n{}^\D$, and recall that $\W$ should satisfy the conditions in (\ref{rhoLRW1}). Then it follows from Lemma \ref{DWcotp} and Equation (\ref{cotp=span}) that $$\W-\C\;\widetilde{\otimes}\;\W-\W\;\widetilde{\otimes}\;\D
\in\bigoplus_{i=1}^{n-1} \bigoplus_{\mathcal{E}\in\mathcal{P}_0}
\bigoplus_{\substack{\X\in{}^\C{\mathcal{P}_i}^\mathcal{E},\\
    \mathcal{Y}\in{}^\mathcal{E}{\mathcal{P}_{n-i}}^\D}}
    \mathrm{span}(\X\;\widetilde{\otimes}\;\Y).$$
However, note that
$$
(\rho^L\otimes\id)(\X\;\widetilde{\otimes}\;\Y)
=\C\;\widetilde{\otimes}\;(\X\;\widetilde{\otimes}\;\Y)\;\;\;\;
\text{and}\;\;\;\;
(\id\otimes\rho^R)(\X\;\widetilde{\otimes}\;\Y)
=(\X\;\widetilde{\otimes}\;\Y)\;\widetilde{\otimes}\;\D,
$$
by applying (\ref{rhoLRW1}) to ${}^\C\mathcal{P}_i{}^\E$ and ${}^\E\mathcal{P}_{n-i}{}^\D$. Thus according to Corollary \ref{NM},
we obtain unique scalars $\alpha_{\mathcal{X},\mathcal{Y}}\in\k$, all but finitely many zero, such that
$$
\Delta(\W)=\C\;\widetilde{\otimes}\;\W+\W\;\widetilde{\otimes}\;\D+\sum_{i=1}^{n-1}
\sum_{\mathcal{E}\in\mathcal{P}_0}\sum_
{\substack{\X\in{}^\C{\mathcal{P}_i}^\mathcal{E},\\
    \mathcal{Y}\in{}^\mathcal{E}{\mathcal{P}_{n-i}}^\D}}
    \alpha_{\X, \Y}(\X\;\widetilde{\otimes}\;\Y).
$$
\end{proof}

\begin{remark}
When $n=1$, Equation (\ref{deltaW}) becomes $\Delta(\W)=\C\;\widetilde{\otimes}\;\W+\W\;\widetilde{\otimes}\;\D$.
Consequently, every $\W\in {}^\C\mathcal{P}_1{}^\D$ is a $(\C, \D)$-primitive matrix in the sense of \cite[Definition 3.2]{LZ19}.
\end{remark}

Finally, the following corollary is provided as a consequence of the coassociativity law for $H$.

\begin{corollary}\label{coassW}
With the notations of Theorem \ref{mainthem}, let $\C, \D\in\mathcal{P}_0$ and $\W\in {}^\C{\mathcal{P}_n}^\D$ where $n\geq1$. Then
\begin{equation}\label{cxyw}
\sum_{i=1}^{n-1}\sum_{\mathcal{E}\in\mathcal{P}_0}
    \sum_{\substack{\mathcal{X}\in{}^\C{\mathcal{P}_i}^\mathcal{E},\\
    \mathcal{Y}\in{}^\mathcal{E}{\mathcal{P}_{n-i}}^\D}}
    \alpha_{\mathcal{X},\mathcal{Y}}\X\;\widetilde{\otimes}\
    \big[\Delta(\Y)-\Y\;\widetilde{\otimes}\;\D\big]
=
\sum_{i=1}^{n-1}\sum_{\mathcal{E}\in\mathcal{P}_0}
    \sum_{\substack{\mathcal{X}\in{}^\C{\mathcal{P}_i}^\mathcal{E},\\
    \mathcal{Y}\in{}^\mathcal{E}{\mathcal{P}_{n-i}}^\D}}
    \alpha_{\mathcal{X},\mathcal{Y}}
    \big[\Delta(\X)-\C\;\widetilde{\otimes}\;\X\big]\;\widetilde{\otimes}\;\Y.
\end{equation}
\end{corollary}

\begin{proof}
%By the coproduct formula for $\W$ given in Theorem \ref{mainthem} (\ref{deltaW}), we know that
%$$
%\Delta(\mathcal{W})=\C\;\widetilde{\otimes}\;\mathcal{W}
%    +\sum_{i=1}^{n-1}\sum_{\mathcal{E}\in\mathcal{P}_0}
%    \sum_{\substack{\mathcal{X}\in{}^\C{\mathcal{P}_i}^\mathcal{E},\\
%    \mathcal{Y}\in{}^\mathcal{E}{\mathcal{P}_{n-i}}^\D}}
%    \alpha_{\mathcal{X},\mathcal{Y}}(\mathcal{X}\;\widetilde{\otimes}\;
%    \mathcal{Y})+\mathcal{W}\;\widetilde{\otimes}\;\D.
%$$
Applying $\id\otimes\Delta$ to Equation $(\ref{deltaW})$ and using $\Delta(\D)=\D\;\widetilde{\otimes}\;\D$,
we obtain
\begin{align*}
(\id\otimes\Delta)\Delta(\W)&=
\C\;\widetilde{\otimes}\;\Delta(\W)+\W\;\widetilde{\otimes}\;\Delta(\D)
+\sum_{i=1}^{n-1}\sum_{\mathcal{E}\in\mathcal{P}_0}
    \sum_{\substack{\mathcal{X}\in{}^\C{\mathcal{P}_i}^\mathcal{E},\\
    \mathcal{Y}\in{}^\mathcal{E}{\mathcal{P}_{n-i}}^\D}}
    \alpha_{\mathcal{X},\mathcal{Y}}(\mathcal{X}\;\widetilde{\otimes}\;
    \Delta(\mathcal{Y})) \\
&=
\C\;\widetilde{\otimes}\;\C\;\widetilde{\otimes}\;\W
+\C\;\widetilde{\otimes}\;\W\;\widetilde{\otimes}\;\D
+\W\;\widetilde{\otimes}\;\D\;\widetilde{\otimes}\;\D\\
&\quad+
\sum_{i=1}^{n-1}\sum_{\mathcal{E}\in\mathcal{P}_0}
\sum_{\substack{\mathcal{X}\in{}^\C{\mathcal{P}_i}^\mathcal{E},\\
    \mathcal{Y}\in{}^\mathcal{E}{\mathcal{P}_{n-i}}^\D}}
    \alpha_{\mathcal{X},\mathcal{Y}}\big[
    \C\;\widetilde{\otimes}\;\X\;\widetilde{\otimes}\;\Y
+\X\;\widetilde{\otimes}\;\Delta(\Y)\big].
\end{align*}
Similarly, applying $\Delta\otimes\id$ to Equation \ref{deltaW} and using $\Delta(\C)=\C\;\widetilde{\otimes}\;\C$, we obtain
\begin{align*}
(\Delta\otimes\id)\Delta(\W)&=
\Delta(\C)\;\widetilde{\otimes}\;\W+\Delta(\W)\;\widetilde{\otimes}\;\D
+\sum_{i=1}^{n-1}\sum_{\mathcal{E}\in\mathcal{P}_0}
    \sum_{\substack{\mathcal{X}\in{}^\C{\mathcal{P}_i}^\mathcal{E},\\
    \mathcal{Y}\in{}^\mathcal{E}{\mathcal{P}_{n-i}}^\D}}
    \alpha_{\mathcal{X},\mathcal{Y}}(\Delta(\mathcal{X})\;\widetilde{\otimes}\;
    \mathcal{Y})
    \\
&=
\C\;\widetilde{\otimes}\;\C\;\widetilde{\otimes}\;\W
+\C\;\widetilde{\otimes}\;\W\;\widetilde{\otimes}\;\D
+\W\;\widetilde{\otimes}\;\D\;\widetilde{\otimes}\;\D\\
&\quad+
\sum_{i=1}^{n-1}\sum_{\mathcal{E}\in\mathcal{P}_0}
\sum_{\substack{\mathcal{X}\in{}^\C{\mathcal{P}_i}^\mathcal{E},\\
    \mathcal{Y}\in{}^\mathcal{E}{\mathcal{P}_{n-i}}^\D}}
    \alpha_{\mathcal{X},\mathcal{Y}}\big[
    \Delta(\X)\;\widetilde{\otimes}\;\Y
+\X\;\widetilde{\otimes}\;\Y\;\widetilde{\otimes}\;\D\big].
\end{align*}
By the coassociativity of $\Delta$, we obtain Equation (\ref{coassW}) as a result.
\end{proof}

\section{Application: Extension formula for the antipode of weak bialgebras}
In this section, we apply the matrix basis decomposition in Theorem \ref{mainthem} to weak bialgebras. Together with Corollary \ref{coassW}, it provides explicit formulas.

We retain the notation of Section 2 and Section 3. In particular, $\k$ is an algebraically closed field. The grading imposed below is a coalgebra grading.

We begin by recalling the definitions of weak bialgebras and their antipodes,
see\cite[Definition 2.1]{BNS99} and \cite[Section 4 and 6]{S03}.

\begin{definition}\label{WHA}
A weak bialgebra  over $\k$ is a quintuple $(H, m, u, \Delta, \varepsilon)$  such that
%satisfying Axioms 1-3 below. If $(A, \delta, u, \Delta, \varepsilon, S)$ satisfies Axioms 1-4 below it is called a weak Hopf algebra $(WHA)$.
$(H, m, u)$ is a  unital associative $\k$-algebra,
    %over a field $\k$ with multiplication $\delta: A\otimes A\rightarrow A$ and unit $u: \k\rightarrow A$. Namely, $\delta$ and $u$ are $\k$-linear and satisfy
%    $$
%    Associativity: \delta\circ(\delta\otimes \id)=\delta\circ(\id\otimes\delta),
%    $$
%    $$
%    Unit\;property: \delta\circ(u\otimes \id)=\id=\delta\circ(\id\otimes u).
%    $$
%Later on we suppress $\delta$ and $u$, just write $xy$ for $\delta(x, y)$, and use the unit element 1:=$u(1)$ instead of $u$.
$(H, \Delta, \varepsilon)$ is a counital coassociative coalgebra, and
 %over $\k$ with comultiplication $\Delta: A\rightarrow A\otimes A$ and counit $\varepsilon: A\rightarrow\k$. Namely, $\Delta$ and $\varepsilon$ are $\k$-linear and satisfy
%    $$
%    Coassocaitivity:(\Delta\otimes\id)\circ\Delta
%    =(\id\otimes\Delta)\circ\Delta,
%    $$
%    $$
%    Counit\;property:(\varepsilon\otimes\id)\circ\Delta
%    =\id=(\id\otimes\varepsilon)\circ\Delta.
%    $$
the following identities hold for all $x, y, z\in H$:
$$\Delta(xy)=\Delta(x)\Delta(y)$$
    $$
    %Weak multiplicativity of the counit:
    \varepsilon(xyz)=\sum\varepsilon(xy_{(1)})
    \varepsilon(y_{(2)}z),\;\;\;\;\varepsilon(xyz)=\sum\varepsilon(xy_{(2)})
    \varepsilon(y_{(1)}z),
    $$
    $$
    %Weak comultiplicativity of the unit:
    \Delta_3(1)=(\Delta(1)\otimes1)(1\otimes\Delta(1)),\;\;\;\;
    \Delta_3(1)=(1\otimes\Delta(1))(\Delta(1)\otimes1).
    $$
%Here $\Delta_3=(\Delta\otimes\id)\circ\Delta=
%(\id\otimes\Delta)\circ\Delta$.
Here $1=u(1_\k)$, $\Delta(1)=\sum1_{(1)}\otimes1_{(2)}$.
The source and target counital maps $\varepsilon_s, \varepsilon_t: H\rightarrow H$ are defined by
\begin{eqnarray*}
\varepsilon_s(x)=\sum 1_{(1)}\varepsilon(x1_{(2)}),\;\;
\varepsilon_t(x)=\sum \varepsilon(1_{(1)}x)1_{(2)}, \;\;\text{for all $x\in H$}.
\end{eqnarray*}

In addition, a weak bialgebra $H$ is a weak Hopf algebra if it admits a linear map $S: H\rightarrow H$ such that, for every $x\in H$,
\begin{equation}\label{idS}
    \sum x_{(1)}S(x_{(2)})=\varepsilon_t(x),\;\;\;\;
    \sum S(x_{(1)})x_{(2)}=\varepsilon_s(x),
    \;\;\;\;
    \sum S(x_{(1)})x_{(2)}S(x_{(3)})=S(x).
\end{equation}
Such a map is called an antipode.
\end{definition}

For a subcoalgebra $K\subseteq H$ and linear maps $f, g: K\rightarrow H$, their convolution product is
$$
(f\ast g)(h)=\sum f(h_{(1)})g(h_{(2)})\;\;\;\;(h\in K).
$$
We shall use the elementary convolution identities
\begin{eqnarray}\label{ttss}
\varepsilon_t\ast\id=\id=\id\ast\varepsilon_s,\;\;
\varepsilon_t\ast\varepsilon_t=\varepsilon_t,\;\;
\varepsilon_s\ast\varepsilon_s=\varepsilon_s.
\end{eqnarray}
The first two identities follow by applying the counit to $\Delta(1)\Delta(h)=\Delta(h)=\Delta(h)\Delta(1)$ for any $h\in H$.
The counital map identities see \cite[Section 4]{S03}
\begin{eqnarray}
\varepsilon_t(x)\varepsilon_t(y)=\varepsilon_t(\varepsilon_t(x)y),
\;\;\;\;\varepsilon_s(x)\varepsilon_s(y)=\varepsilon_s(x\varepsilon_s(y)),
\;\;\;\;(\forall x, y\in H)
\end{eqnarray}
then imply
$$
(\varepsilon_t\ast\varepsilon_t)(h)=\sum\varepsilon_t(h_{(1)})
\varepsilon_t(h_{(2)})=\varepsilon_t
(\sum\varepsilon_t(h_{(1)})h_{(2)})=\varepsilon_t(h),
$$
$$
(\varepsilon_s\ast\varepsilon_s)(h)=\sum\varepsilon_s(h_{(1)})
\varepsilon_s(h_{(2)})=\varepsilon_s
(\sum h_{(1)}\varepsilon_s(h_{(2)}))=\varepsilon_s(h).
$$

For the remainder of this section, let $H=\bigoplus_{n\geq0}H(n)$ be a weak bialgebra. As a coalgebra, $H$ is $\mathbb{N}$-graded with $H(0)=H_0$, as in Section 3.2.

The following lemma follows immediately from Corollary \ref{coassW}.
\begin{lemma}\label{identitys1}
With the notation of Theorem \ref{mainthem}.
%Let $H=\bigoplus_{n\geq0}H(n)$ be a weak bialgebra whose underlying coalgebra is graded.
Suppose that $S: \bigoplus_{i=0}^{n-1}H(i)\rightarrow H$ is  a linear map satisfying $\id\ast S=\varepsilon_t$ and $S\ast\id=\varepsilon_s$ on $\bigoplus_{i=0}^{n-1}H(i)$, then
\begin{equation}\label{xycd}
\sum_{i=1}^{n-1}\sum_{\mathcal{E}\in\mathcal{P}_0}
\sum_{\substack{\mathcal{X}\in{}^\C\mathcal{P}_i^\mathcal{E},\\ \mathcal{Y}\in{}^\mathcal{E}\mathcal{P}_{n-i}^\D}}
\alpha_{\mathcal{X},\mathcal{Y}}
S(\mathcal{X})\big[\varepsilon_t (\mathcal{Y})-\mathcal{Y}S(\D)\big]
= \sum_{i=1}^{n-1}\sum_{\mathcal{E}\in\mathcal{P}_0}\sum_{\substack{\mathcal{X}
\in{}^\C\mathcal{P}_i^\mathcal{E},\\ \mathcal{Y}\in{}^\mathcal{E}\mathcal{P}_{n-i}^\D}}
\alpha_{\mathcal{X},\mathcal{Y}}\big[
\varepsilon_s(\mathcal{X})-S(\C)\mathcal{X}
\big]S(\mathcal{Y}).
\end{equation}
\end{lemma}

\begin{proof}
Applying $m\circ(m\otimes\id)\circ(S\otimes\id\otimes S)$ to both sides of the Equation (\ref{cxyw}) and using (\ref{idS}),
we obtain Equation (\ref{xycd}).
\end{proof}

In addition, under the hypotheses of Lemma \ref{identitys1} , if $S\ast\id\ast S=S$ on $\bigoplus_{i=0}^{n-1}H(i)$, then
\begin{eqnarray}\label{sSSt}
\varepsilon_s\ast S=S=S\ast \varepsilon_t \;\;\;\; \text{on $\bigoplus_{i=0}^{n-1}H(i)$}.
\end{eqnarray}

\begin{corollary}\label{CXYXYD}
Under the hypotheses of Lemma \ref{identitys1}, assume in addition that $S\ast\id\ast S= S$ on $\bigoplus_{i=0}^{n-1}H(i)$. Then
%Let $H=\bigoplus_{n\geq0}H(n)$ be a coradically graded coalgebra and $\mathcal{W}\in{}^\C\mathcal{P}_n^\D$. Then the following identity holds:
\begin{equation}\label{CXY=XYD}
\sum_{i=1}^{n-1}\sum_{\mathcal{E}\in\mathcal{P}_0}
    \sum_{\substack{\mathcal{X}\in{}^\C{\mathcal{P}_i}^\mathcal{E},\\
    \mathcal{Y}\in{}^\mathcal{E}{\mathcal{P}_{n-i}}^\D}}
    \alpha_{\mathcal{X},\mathcal{Y}}
    \X\varepsilon_s(\Y)S(\D)
~=~
\sum_{i=1}^{n-1}\sum_{\mathcal{E}\in\mathcal{P}_0}
    \sum_{\substack{\mathcal{X}\in{}^\C{\mathcal{P}_i}^\mathcal{E},\\
    \mathcal{Y}\in{}^\mathcal{E}{\mathcal{P}_{n-i}}^\D}}
    \alpha_{\mathcal{X},\mathcal{Y}}
\C\varepsilon_s(\X)S(\Y)
\end{equation}
and
\begin{equation}\label{XYD}
\sum_{i=1}^{n-1}\sum_{\mathcal{E}\in\mathcal{P}_0}
    \sum_{\substack{\mathcal{X}\in{}^\C{\mathcal{P}_i}^\mathcal{E},\\
    \mathcal{Y}\in{}^\mathcal{E}{\mathcal{P}_{n-i}}^\D}}
    \alpha_{\mathcal{X},\mathcal{Y}}
\varepsilon_t(\X)[\varepsilon_t(\Y)-
\Y S(\D)]
=
\sum_{i=1}^{n-1}\sum_{\mathcal{E}\in\mathcal{P}_0}
    \sum_{\substack{\mathcal{X}\in{}^\C{\mathcal{P}_i}^\mathcal{E},\\
    \mathcal{Y}\in{}^\mathcal{E}{\mathcal{P}_{n-i}}^\D}}
    \alpha_{\mathcal{X},\mathcal{Y}}
    [\X-\varepsilon_t(\C)\X]S(\Y).
\end{equation}
\end{corollary}

\begin{proof}
\begin{itemize}
\item[(1)]
Applying $m_3\circ(\id\otimes\;\varepsilon_s\otimes S)$ to both sides of (\ref{cxyw}) and using (\ref{ttss}), (\ref{sSSt}), we obtain
\begin{eqnarray*}
\sum_{i=1}^{n-1}\sum_{\mathcal{E}\in\mathcal{P}_0}
    \sum_{\substack{\mathcal{X}\in{}^\C{\mathcal{P}_i}^\mathcal{E},\\
    \mathcal{Y}\in{}^\mathcal{E}{\mathcal{P}_{n-i}}^\D}}
    \alpha_{\mathcal{X},\mathcal{Y}}
    \X[S(\Y)-\varepsilon_s(\Y)S(\D)]
=
\sum_{i=1}^{n-1}\sum_{\mathcal{E}\in\mathcal{P}_0}
    \sum_{\substack{\mathcal{X}\in{}^\C{\mathcal{P}_i}^\mathcal{E},\\
    \mathcal{Y}\in{}^\mathcal{E}{\mathcal{P}_{n-i}}^\D}}
    \alpha_{\mathcal{X},\mathcal{Y}}
    \big[\X-\C\varepsilon_s(\X)\big]S(\Y).
%&=&
%\sum_{i=1}^{n-1}\sum_{\mathcal{E}\in\mathcal{P}_0}
%    \sum_{\substack{\mathcal{X}\in{}^\C{\mathcal{P}_i}^\mathcal{E},\\
%    \mathcal{Y}\in{}^\mathcal{E}{\mathcal{P}_{n-i}}^\D}}
%    \alpha_{\mathcal{X},\mathcal{Y}}
%\big[\X S(\Y)-\X S(\Y)\big]=0.
\end{eqnarray*}
This gives (\ref{CXY=XYD}).
\item[(2)] Applying $m_3\circ(\varepsilon_t\otimes
\id\otimes S)$ directly to Equation (\ref{cxyw}) and using (\ref{idS}), (\ref{ttss}), we directly obtain (\ref{XYD}).
%\begin{eqnarray*}
%&&
%\sum_{i=1}^{n-1}\sum_{\mathcal{E}\in\mathcal{P}_0}
%    \sum_{\substack{\mathcal{X}\in{}^\C{\mathcal{P}_i}^\mathcal{E},\\
%    \mathcal{Y}\in{}^\mathcal{E}{\mathcal{P}_{n-i}}^\D}}
%    \alpha_{\mathcal{X},\mathcal{Y}}
%[\varepsilon_t(\X)\Y S(\D)-\varepsilon_t(\C)\X S(\Y)] \nonumber \\
%&\overset{}{=}&
%\sum_{i=1}^{n-1}\sum_{\mathcal{E}\in\mathcal{P}_0}
%    \sum_{\substack{\mathcal{X}\in{}^\C{\mathcal{P}_i}^\mathcal{E},\\
%    \mathcal{Y}\in{}^\mathcal{E}{\mathcal{P}_{n-i}}^\D}}
%    \alpha_{\mathcal{X},\mathcal{Y}}
%    [\varepsilon_t(\X)\varepsilon_t(\Y)-\X S(\Y)].
%\end{eqnarray*}
\end{itemize}
\end{proof}

\begin{remark}\label{SXYD}
Under the same hypotheses, applying $m_3\circ(S\otimes \varepsilon_t\otimes
\id$) and $m_3\circ(S\otimes \id\otimes \varepsilon
_s)$ to (\ref{cxyw}) gives, respectively,
\begin{equation}\label{XYD=CXY}
\sum_{i=1}^{n-1}\sum_{\mathcal{E}\in\mathcal{P}_0}
    \sum_{\substack{\mathcal{X}\in{}^\C{\mathcal{P}_i}^\mathcal{E},\\
    \mathcal{Y}\in{}^\mathcal{E}{\mathcal{P}_{n-i}}^\D}}
    \alpha_{\mathcal{X},\mathcal{Y}}
    S(\C)\varepsilon_t(\X)\Y
~=~
\sum_{i=1}^{n-1}\sum_{\mathcal{E}\in\mathcal{P}_0}
    \sum_{\substack{\mathcal{X}\in{}^\C{\mathcal{P}_i}^\mathcal{E},\\
    \mathcal{Y}\in{}^\mathcal{E}{\mathcal{P}_{n-i}}^\D}}
    \alpha_{\mathcal{X},\mathcal{Y}}
S(\X)\varepsilon_t(\Y)\D.
\end{equation}
and
\begin{equation}\label{CXY}
\sum_{i=1}^{n-1}\sum_{\mathcal{E}\in\mathcal{P}_0}
    \sum_{\substack{\mathcal{X}\in{}^\C{\mathcal{P}_i}^\mathcal{E},\\
    \mathcal{Y}\in{}^\mathcal{E}{\mathcal{P}_{n-i}}^\D}}
    \alpha_{\mathcal{X},\mathcal{Y}}
\big[S(\C)\X-\varepsilon_s(\X)\big]\varepsilon_s(\Y)
=
\sum_{i=1}^{n-1}\sum_{\mathcal{E}\in\mathcal{P}_0}
    \sum_{\substack{\mathcal{X}\in{}^\C{\mathcal{P}_i}^\mathcal{E},\\
    \mathcal{Y}\in{}^\mathcal{E}{\mathcal{P}_{n-i}}^\D}}
    \alpha_{\mathcal{X},\mathcal{Y}}
    S(\X)\big[\Y\varepsilon_s(\D)-\Y\big].
\end{equation}
\end{remark}

Moreover,
applying $\varepsilon_t\ast\id=\id$ and $\id\ast\varepsilon_s=\id$ to (\ref{deltaW}) gives, respectively,
\begin{eqnarray*}
\W-\varepsilon_t(\C)\W=\varepsilon_t(\W)\D+
\sum_{i=1}^{n-1}\sum_{\mathcal{E}\in\mathcal{P}_0}
    \sum_{\substack{\mathcal{X}\in{}^\C{\mathcal{P}_i}^\mathcal{E},\\
    \mathcal{Y}\in{}^\mathcal{E}{\mathcal{P}_{n-i}}^\D}}
    \alpha_{\mathcal{X},\mathcal{Y}}
    \varepsilon_t(\X)\Y, \\
\W-\W\varepsilon_s(\D)=\C\varepsilon_s(\W)+
\sum_{i=1}^{n-1}\sum_{\mathcal{E}\in\mathcal{P}_0}
    \sum_{\substack{\mathcal{X}\in{}^\C{\mathcal{P}_i}^\mathcal{E},\\
    \mathcal{Y}\in{}^\mathcal{E}{\mathcal{P}_{n-i}}^\D}}
    \alpha_{\mathcal{X},\mathcal{Y}}
    \X\varepsilon_s(\Y).
\end{eqnarray*}
Subtracting the two identities, we obtain
\begin{equation}\label{es-et}
\C\varepsilon_s(\W)-\varepsilon_t(\C)\W=\varepsilon_t(\W)\D
-\W\varepsilon_s(\D)+\sum_{i=1}^{n-1}\sum_{\mathcal{E}\in\mathcal{P}_0}
    \sum_{\substack{\mathcal{X}\in{}^\C{\mathcal{P}_i}^\mathcal{E},\\
    \mathcal{Y}\in{}^\mathcal{E}{\mathcal{P}_{n-i}}^\D}}
    \alpha_{\mathcal{X},\mathcal{Y}}(\varepsilon_t(\X)\Y-\X\varepsilon_s(\Y)).\upsilon
\end{equation}

We now use these identities to extend $S$ from $\bigoplus_{i=0}^{n-1}H(i)$ to $H$.
Lemma \ref{identitys1} shows that the following two formulas for $S$ agree. Corollary \ref{CXYXYD} and Remark \ref{SXYD} provide the identities needed to verify the antipode axioms.

\begin{proposition}\label{Sn}
%Fix $n\geq1$, and suppose that $S: \bigoplus_{i=0}^{n-1}H(i)\rightarrow H$ is a linear map satisfying
Under the hypotheses of Corollary \ref{CXYXYD},
%$$
%\id\ast S=\varepsilon_t,\;\;S\ast\id=\varepsilon_s,\;\;
%S\ast\id\ast S=S\;\;\;\;\text{on}\;\; \bigoplus_{i=0}^{n-1}H(i).
%$$
for each $\C, \D\in\mathcal{P}_0$ and $\W\in{}^\C\mathcal{P}_n{}^\D$,
define the values of $S$ on the entries of $\W$ by
\begin{align}\label{antipode1}
S(\W)&=S(\C)\varepsilon_t(\W)+\varepsilon_s(\W)S(\D)-S(\C)\W S(\D)  \nonumber\\
&\quad +
\sum_{i=1}^{n-1}\sum_{\mathcal{E}\in\mathcal{P}_0}
    \sum_{\substack{\mathcal{X}\in{}^\C{\mathcal{P}_i}^\mathcal{E},\\
    \mathcal{Y}\in{}^\mathcal{E}{\mathcal{P}_{n-i}}^\D}}
    \alpha_{\mathcal{X},\mathcal{Y}}
    \big[\varepsilon_s(\X)S(\Y)-S(\C)\X S(\Y)\big].
\end{align}
or equivalently,
\begin{align}\label{antipode2}
S(\W)&=S(\C)\varepsilon_t(\W)+\varepsilon_s(\W)S(\D)-S(\C)\W S(\D)  \nonumber\\
&\quad +
\sum_{i=1}^{n-1}\sum_{\mathcal{E}\in\mathcal{P}_0}
\sum_{\substack{\mathcal{X}\in{}^\C\mathcal{P}_i^\mathcal{E},\\ \mathcal{Y}\in{}^\mathcal{E}\mathcal{P}_{n-i}^\D}}
\alpha_{\mathcal{X},\mathcal{Y}}\big[ S(\mathcal{X})\varepsilon_t (\mathcal{Y})-S(\mathcal{X})\mathcal{Y}S(\D)\big].
\end{align}
The given map extends uniquely to a linear map $S: H\rightarrow H$, which satisfies
$$
\id\ast S=\varepsilon_t,\;\;S\ast\id=\varepsilon_s,\;\;
S\ast\id\ast S=S\;\;\;\;\text{on $H$}.
$$
\end{proposition}

\begin{proof}
By Theorem \ref{mainthem}, the entries of all matrices in $\mathcal{P}_n$ form a basis of $H(n)$. Since $H=(\bigoplus_{i=0}^{n-1}H(i))\oplus H(n)$, the formula (\ref{antipode1}) therefore defines a unique linear extension. We continue to denote it by $S$.

Fix $\C, \D\in \mathcal{P}_0$ and $\W\in {}^\C\mathcal{P}_n{}^\D$.
To verify $\id\ast S=\varepsilon_t$, we use (\ref{deltaW}) to obtain
%(\ref{antipode1})
\begin{align*}
&
\C S(\W)+\W S(\D)+\sum_{i=1}^{n-1}\sum_{\mathcal{E}\in\mathcal{P}_0}
    \sum_{\substack{\mathcal{X}\in{}^\C{\mathcal{P}_i}^\mathcal{E},\\
    \mathcal{Y}\in{}^\mathcal{E}{\mathcal{P}_{n-i}}^\D}}
    \alpha_{\mathcal{X},\mathcal{Y}}\X S(\Y) \\
%&=
%\C S(\C)\varepsilon_t(\W)+\C\varepsilon_s(\W)S(\D)-\C S(\C)\W S(\D)+\W S(\D) \\
%&\quad+\sum_{i=1}^{n-1}\sum_{\mathcal{E}\in\mathcal{P}_0}
%    \sum_{\substack{\mathcal{X}\in{}^\C{\mathcal{P}_i}^\mathcal{E},\\
%    \mathcal{Y}\in{}^\mathcal{E}{\mathcal{P}_{n-i}}^\D}}
%    \alpha_{\mathcal{X},\mathcal{Y}}
%    \big[\C\varepsilon_s(\X)S(\Y)-\C S(\C)\X S(\Y)+\X S(\Y)
%    \big] \\
%&=
%\varepsilon_t(\C)\varepsilon_t(\W)+\C\varepsilon_s(\W)S(\D)-
%\varepsilon_t(\C)\W S(\D)+\W S(\D) \\
%&\quad+\sum_{i=1}^n\sum_{\mathcal{E}\in\mathcal{P}_0}
%    \sum_{\substack{\mathcal{X}\in{}^\C{\mathcal{P}_i}^\mathcal{E},\\
%    \mathcal{Y}\in{}^\mathcal{E}{\mathcal{P}_{n-i}}^\D}}
%    \alpha_{\mathcal{X},\mathcal{Y}}
%    \big[\C\varepsilon_s(\X)S(\Y)-\varepsilon_t(\C)\X S(Y)+\X S(Y)
%    \big] \\
&\overset{(\ref{antipode1}),(\ref{idS})}{=}
\varepsilon_t(\C)\varepsilon_t(\W)+[\C\varepsilon_s(\W)
+\W-\varepsilon_t(\C)\W ]S(\D) \\
&\quad+\sum_{i=1}^{n-1}\sum_{\mathcal{E}\in\mathcal{P}_0}
    \sum_{\substack{\mathcal{X}\in{}^\C{\mathcal{P}_i}^\mathcal{E},\\
    \mathcal{Y}\in{}^\mathcal{E}{\mathcal{P}_{n-i}}^\D}}
    \alpha_{\mathcal{X},\mathcal{Y}}
    \big[\C\varepsilon_s(\X)S(\Y)-\varepsilon_t(\C)\X S(\Y)+\X S(\Y)
    \big] \\
%&=
%\varepsilon_t(\C)\varepsilon_t(\W)+\C\varepsilon_s(\W)S(\D)
%+\varepsilon_t(\W)\varepsilon_t(\D) \\
%&\quad+\sum_{i=1}^{n-1}\sum_{\mathcal{E}\in\mathcal{P}_0}
%    \sum_{\substack{\mathcal{X}\in{}^\C{\mathcal{P}_i}^\mathcal{E},\\
%    \mathcal{Y}\in{}^\mathcal{E}{\mathcal{P}_{n-i}}^\D}}
%    \alpha_{\mathcal{X},\mathcal{Y}}
%    \big[\C\varepsilon_s(\X)S(\Y)-\varepsilon_t(\C)\X S(\Y)+\X S(\Y)+\varepsilon_t(\X)\Y S(\D)
%    \big] \\
&\overset{(\ref{es-et}),(\ref{sSSt})}=
\varepsilon_t(\C)\varepsilon_t(\W)+\varepsilon_t(\W)\varepsilon_t(\D)
%+\W S(\D)-\W\varepsilon_s(\D) S(\D)
+\sum_{i=1}^{n-1}\sum_{\mathcal{E}\in\mathcal{P}_0}
    \sum_{\substack{\mathcal{X}\in{}^\C{\mathcal{P}_i}^\mathcal{E},\\
    \mathcal{Y}\in{}^\mathcal{E}{\mathcal{P}_{n-i}}^\D}}
    \alpha_{\mathcal{X},\mathcal{Y}}\\
&\quad
    \big[\C\varepsilon_s(\X)S(\Y)-\X\varepsilon_s(\Y)S(\D)
    +\varepsilon_t(\X)\Y S(\D)-\varepsilon_t(\C)\X S(\Y)+\X S(\Y)\big] \\
&\overset{(\ref{CXY=XYD}),(\ref{XYD})}{=}
\varepsilon_t(\C)\varepsilon_t(\W)+\varepsilon_t(\W)
\varepsilon_t(\D)
+\sum_{i=1}^{n-1}\sum_{\mathcal{E}\in\mathcal{P}_0}
\sum_{\substack{\mathcal{X}\in{}^\C{\mathcal{P}_i}^\mathcal{E},\\
\mathcal{Y}\in{}^\mathcal{E}{\mathcal{P}_{n-i}}^\D}}
\alpha_{\mathcal{X},\mathcal{Y}}
\varepsilon_t(\X)\varepsilon_t(\Y) \\
&= \varepsilon_t(\W).
\end{align*}
Next, using the equations (\ref{deltaW}) and (\ref{antipode2}), we obtain
%\begin{align*}
%S(\W)
%&=S(\C)\varepsilon_t(\W)+\varepsilon_s(\W)S(\D)-S(\C)\W S(\D)  \nonumber \\
%&\quad +
%\sum_{i=1}^{n-1}\sum_{\mathcal{E}\in\mathcal{P}_0}
%\sum_{\substack{\mathcal{X}\in{}^\C\mathcal{P}_i^\mathcal{E},\\ \mathcal{Y}\in{}^\mathcal{E}\mathcal{P}_{n-i}^\D}}
%\alpha_{\mathcal{X},\mathcal{Y}}\big[ S(\mathcal{X})\varepsilon_t (\mathcal{Y})-S(\mathcal{X})\mathcal{Y}S(\D)\big],
%\end{align*}
\begin{align*}
&
S(\C)\W+S(\W)\D
+\sum_{i=1}^{n-1}\sum_{\mathcal{E}\in\mathcal{P}_0}
    \sum_{\substack{\mathcal{X}\in{}^\C{\mathcal{P}_i}^\mathcal{E},\\
    \mathcal{Y}\in{}^\mathcal{E}{\mathcal{P}_{n-i}}^\D}}
    \alpha_{\mathcal{X},\mathcal{Y}}S(\X)\Y \\
%&=
%S(\C)\W+S(\C)\varepsilon_t(\W)\D+\varepsilon_s(\W)S(\D)\D
%-S(\C)\W S(\D)\D\\
%&\quad+
%\sum_{i=1}^n\sum_{\mathcal{E}\in\mathcal{P}_0}
%\sum_{\substack{\mathcal{X}\in{}^\C\mathcal{P}_i^\mathcal{E},\\ \mathcal{Y}\in{}^\mathcal{E}\mathcal{P}_{n-i}^\D}}
%\alpha_{\mathcal{X},\mathcal{Y}}\big[ S(\mathcal{X})\varepsilon_t (\mathcal{Y})\D-S(\mathcal{X})\mathcal{Y}S(\D)\D+S(\X)\Y\big] \\
&\overset{(\ref{antipode2}),(\ref{idS})}{=}
S(\C)[\W-\W \varepsilon_s(\D)+\varepsilon_t(\W)\D]
+\varepsilon_s(\W)\varepsilon_s(\D)\\
&\quad+
\sum_{i=1}^{n-1}\sum_{\mathcal{E}\in\mathcal{P}_0}
\sum_{\substack{\mathcal{X}\in{}^\C\mathcal{P}_i^\mathcal{E},\\ \mathcal{Y}\in{}^\mathcal{E}\mathcal{P}_{n-i}^\D}}
\alpha_{\mathcal{X},\mathcal{Y}}\big[ S(\mathcal{X})\varepsilon_t (\mathcal{Y})\D-S(\mathcal{X})\mathcal{Y}\varepsilon_s(\D)+S(\X)\Y
\big] \\
%&=
%S(\C)\C\varepsilon_s(\D)+S(\C)\varepsilon_t(\W)\D
%+\varepsilon_s(\W)\varepsilon_s(\D)\\
%&\quad+
%\sum_{i=1}^n\sum_{\mathcal{E}\in\mathcal{P}_0}
%\sum_{\substack{\mathcal{X}\in{}^\C\mathcal{P}_i^\mathcal{E},\\ \mathcal{Y}\in{}^\mathcal{E}\mathcal{P}_{n-i}^\D}}
%\alpha_{\mathcal{X},\mathcal{Y}}\big[ S(\mathcal{X})\varepsilon_t (\mathcal{Y})\D-S(\mathcal{X})\mathcal{Y}\varepsilon_s(\D)+S(\X)\Y
%+S(\C)\X\varepsilon_s(\Y)\big] \\
%&=
%\varepsilon_s   (\C)\varepsilon_s(\W)+S(\C)\varepsilon_t(\W)\D
%+\varepsilon_s(\W)\varepsilon_s(\D)\\
%&\quad+
%\sum_{i=1}^{n-1}\sum_{\mathcal{E}\in\mathcal{P}_0}
%\sum_{\substack{\mathcal{X}\in{}^\C\mathcal{P}_i^\mathcal{E},\\ \mathcal{Y}\in{}^\mathcal{E}\mathcal{P}_{n-i}^\D}}
%\alpha_{\mathcal{X},\mathcal{Y}}\big[ S(\mathcal{X})\varepsilon_t (\mathcal{Y})\D-S(\mathcal{X})\mathcal{Y}\varepsilon_s(\D)+S(\X)\Y
%+S(\C)\X\varepsilon_s(\Y)\big] \\
&\overset{(\ref{es-et}),(\ref{sSSt})}{=}
\varepsilon_s(\C)\varepsilon_s(\W)
+\varepsilon_s(\W)\varepsilon_s(\D)+
\sum_{i=1}^{n-1}\sum_{\mathcal{E}\in\mathcal{P}_0}
\sum_{\substack{\mathcal{X}\in{}^\C\mathcal{P}_i^\mathcal{E},\\ \mathcal{Y}\in{}^\mathcal{E}\mathcal{P}_{n-i}^\D}}
\alpha_{\mathcal{X},\mathcal{Y}}\\
&\quad+
\big[ S(\mathcal{X})\varepsilon_t (\mathcal{Y})\D-S(\mathcal{C})\varepsilon_t(\X)\Y
+S(\C)\X\varepsilon_s(\Y)-S(\X)\Y\varepsilon_s(\D)+S(\X)\Y\big]\\
&\overset{(\ref{SXYD}),(\ref{XYD=CXY})}{=}
S(\C)\C\varepsilon_s(\W)
+\varepsilon_s(\W)\varepsilon_s(\D)+
\sum_{i=1}^{n-1}\sum_{\mathcal{E}\in\mathcal{P}_0}
\sum_{\substack{\mathcal{X}\in{}^\C\mathcal{P}_i^\mathcal{E},\\ \mathcal{Y}\in{}^\mathcal{E}\mathcal{P}_{n-i}^\D}}
\alpha_{\mathcal{X},\mathcal{Y}}
\varepsilon_s(\X)\varepsilon_s(\Y) \\
&=\varepsilon_s(\W)
\end{align*}
which verifies $S\ast\id=\varepsilon_s$.

Finally,
%we verify the third antipode identity. Applying $\Delta\otimes\id$ to (\ref{deltaW}) and expanding $\Delta(\W)$ gives
%\begin{align*}
%&
%\Delta_3(\W)=\C\;\widetilde{\otimes}\;\C\;\widetilde{\otimes}
%\;\W+\C\;\widetilde{\otimes}\;\W\;\widetilde{\otimes}\;\D
%+\W\;\widetilde{\otimes}\;\D\;\widetilde{\otimes}\;\D \\
%&\quad+
%\sum_{i=1}^{n-1}\sum_{\mathcal{E}\in\mathcal{P}_0}
%\sum_{\substack{\mathcal{X}\in{}^\C\mathcal{P}_i^\mathcal{E},\\ \mathcal{Y}\in{}^\mathcal{E}\mathcal{P}_{n-i}^\D}}
%\alpha_{\mathcal{X},\mathcal{Y}}
%[\Delta(\X)\;\widetilde{\otimes}\;\Y+\X\;\widetilde{\otimes}\;
%\Y\;\widetilde{\otimes}\;\D].
%\end{align*}
applying $S\ast\id\ast S$ to $\W$, and using $S\ast\id=\varepsilon_s$, we obtain
\begin{align*}
&
(S\ast\id\ast S)(\W)=
S(\C)\C S(\W)+S(\C)\W S(\D)
+S(\W)\D S(\D) \\
&\quad+
\sum_{i=1}^{n-1}\sum_{\mathcal{E}\in\mathcal{P}_0}
\sum_{\substack{\mathcal{X}\in{}^\C\mathcal{P}_i^\mathcal{E},\\ \mathcal{Y}\in{}^\mathcal{E}\mathcal{P}_{n-i}^\D}}
\alpha_{\mathcal{X},\mathcal{Y}}
[\varepsilon_s(\X)S(\Y)+S(\X)\Y S(\D)].
\end{align*}
%The left hand
%\begin{align*}
%&
%S(\C)\C S(\W)+S(\C)\W S(\D)
%+S(\W)\D S(\D) \\
%&\quad+
%\sum_{i=1}^n\sum_{\mathcal{E}\in\mathcal{P}_0}
%\sum_{\substack{\mathcal{X}\in{}^\C\mathcal{P}_i^\mathcal{E},\\ \mathcal{Y}\in{}^\mathcal{E}\mathcal{P}_{n-i}^\D}}
%\alpha_{\mathcal{X},\mathcal{Y}}
%[\varepsilon_s(\X)\Y+S(\X)\Y S(\D)]\\
%&=
%\varepsilon_s(\C)S(\W)+S(\C)\W S(\D)+S(\W)\varepsilon_t(\D)\\
%&\quad+
%\sum_{i=1}^n\sum_{\mathcal{E}\in\mathcal{P}_0}
%\sum_{\substack{\mathcal{X}\in{}^\C\mathcal{P}_i^\mathcal{E},\\ \mathcal{Y}\in{}^\mathcal{E}\mathcal{P}_{n-i}^\D}}
%\alpha_{\mathcal{X},\mathcal{Y}}
%[\varepsilon_s(\X)\Y+S(\X)\Y S(\D)]\\
%\end{align*}
%By $\id\ast S=\varepsilon_t$ and $S\ast\id=\varepsilon_s$, which have already been established on $H_n$, we have
%\begin{eqnarray}\label{etW}
%\varepsilon_t(\W)=\C S(\W)+\W S(\D)+
%\sum_{i=1}^{n-1}\sum_{\mathcal{E}\in\mathcal{P}_0}
%\sum_{\substack{\mathcal{X}\in{}^\C\mathcal{P}_i^\mathcal{E},\\ \mathcal{Y}\in{}^\mathcal{E}\mathcal{P}_{n-i}^\D}}
%\alpha_{\mathcal{X},\mathcal{Y}}
%\X S(\Y)
%\end{eqnarray}
%and
%\begin{eqnarray}\label{esW}
%\varepsilon_s(\W)=S(\C)\W+S(\W)\D+
%\sum_{i=1}^{n-1}\sum_{\mathcal{E}\in\mathcal{P}_0}
%\sum_{\substack{\mathcal{X}\in{}^\C\mathcal{P}_i^\mathcal{E},\\ \mathcal{Y}\in{}^\mathcal{E}\mathcal{P}_{n-i}^\D}}
%\alpha_{\mathcal{X},\mathcal{Y}}
%S(\X)\Y.
%\end{eqnarray}
Hence, substituting the previous equations $\varepsilon_t(\W)$ and $\varepsilon_s(\W)$ to (\ref{antipode1}), we obtain $S=S\ast\id\ast S$. Namely,
\begin{align*}
&
S(\W)=
S(\C)\varepsilon_t(\W)+\varepsilon_s(\W)S(\D)-S(\C)\W S(\D)\\
&\quad+
\sum_{i=1}^{n-1}\sum_{\mathcal{E}\in\mathcal{P}_0}
\sum_{\substack{\mathcal{X}\in{}^\C\mathcal{P}_i^\mathcal{E},\\ \mathcal{Y}\in{}^\mathcal{E}\mathcal{P}_{n-i}^\D}}
\alpha_{\mathcal{X},\mathcal{Y}}
[\varepsilon_s(\X)S(\Y)-S(\C)\X S(\Y)] \\
%&=
%S(\C)\C S(\W)+S(\C)\W S(\D)+S(\W)\D S(\D) \\
%&\quad+
%\sum_{i=1}^{n-1}\sum_{\mathcal{E}\in\mathcal{P}_0}
%\sum_{\substack{\mathcal{X}\in{}^\C\mathcal{P}_i^\mathcal{E},\\ \mathcal{Y}\in{}^\mathcal{E}\mathcal{P}_{n-i}^\D}}
%\alpha_{\mathcal{X},\mathcal{Y}}
%[S(\X)\Y S(\D)+\varepsilon_s(\X)S(\Y)] \\
&=
\varepsilon_s(\C)S(\W)+S(\C)\W S(\D)+S(\W)\varepsilon_t(\D)\\
&\quad+
\sum_{i=1}^{n-1}\sum_{\mathcal{E}\in\mathcal{P}_0}
\sum_{\substack{\mathcal{X}\in{}^\C\mathcal{P}_i^\mathcal{E},\\ \mathcal{Y}\in{}^\mathcal{E}\mathcal{P}_{n-i}^\D}}
\alpha_{\mathcal{X},\mathcal{Y}}
[\varepsilon_s(\X)S(\Y)+S(\X)\Y S(\D)]. \\
&=
(S\ast\id\ast S)(\W).
\end{align*}
Thus all three antipode identities hold on $H(n)$.
\end{proof}

\begin{theorem}\label{WHA MT}
Under the same hypotheses.
if there exists a linear map
$S: H(0)\rightarrow H$ such that
$$
S\ast\id=\varepsilon_s,\;\;\id\ast S=\varepsilon_t,\;\;
S\ast\id\ast S=S,
$$
then $H$ is a weak Hopf algebra.
\end{theorem}

\begin{proof}
The given map satisfies all three antipode identities on $H(0)$. Applying Lemma \ref{Sn} successively, we obtain linear extensions to $\bigoplus_{i=0}^nH(i)$ for all $n\geq0$, each satisfying the three identities. Since $H=\bigoplus_{i=0}^nH(i)$, these extensions determine a linear map $S: H\rightarrow H$. Hence the three identities hold on $H$. Thus $S$ is an antipode, and $H$ is a weak Hopf algebra.
\end{proof}

\begin{remark}
This result was mentioned in an unpublished work by H. Huang, X. Wang and J.J. Zhang in the situation when $H_0$ is a weak Hopf subalgebra of $H$.
However, our method includes recursive formulas (\ref{antipode1}) and (\ref{antipode1}) for the antipode, which may be helpful to calculate in specific.
\end{remark}

\section*{Acknowledgement}

The authors are grateful to Professors Gongxiang Liu and Hongdi Huang for useful backgrounds and suggestions.

%AI tools were only used for grammar checking and language editing. The author takes sole responsibility for the contents of the
%paper.

%\section*{References}

\end{document}